# Homoclinic tangencies and hyperbolicity for surface diffeomorphisms

By Enrique R. Pujals and Martín Sambarino*

**Abstract**

We prove here that in the complement of the closure of the hyperbolic surface diffeomorphisms, the ones exhibiting a homoclinic tangency are $C^1$ dense. This represents a step towards the global understanding of dynamics of surface diffeomorphisms.

## 1. Introduction

A long-time goal in the theory of dynamical systems is to describe the dynamics of "big sets" (generic or residual, dense, etc.) in the space of all dynamical systems. It was thought in the sixties that this could be realized by the stable ones, say on compact smooth manifolds without boundary $M$: the set of the structurally stable dynamical systems would form an open and dense subset of all dynamics, endowed with the $C^r$ topology, for some $r \geq 1$, where a system is said to be structurally stable if all nearby ones are conjugate to it through a homeomorphism of the ambient manifold. Nevertheless, examples by Smale [S1] and others showed that this was in fact impossible: there are open sets in the space of diffeomorphisms formed by systems which are not structurally stable. In fact, even $\Omega$-stable systems are not dense ([AS]). Here $\Omega$ stands for the nonwandering set of the dynamics and $\Omega$-stability for stability restricted to the nonwandering set.

The important case of surface diffeomorphisms was not included in such counter-examples. By the end of the sixties, even $\Omega$-stable diffeomorphisms on surfaces were shown not to be dense in the $C^r$ topology with $r \geq 2$ (the case $r = 1$ is still an open problem). This has been done through very original examples of Newhouse (see [N1], [N2], [N3]). In fact, he proved that the unfolding of a homoclinic tangency leads to very rich dynamics: residual subsets

*The first author was partially supported by UFRJ/FAPERJ-Brazil and CNPq/PRONEX-Dyn. Syst. The second author was partially supported by IMPA/CNPq-Brazil, CNPq/PRONEX Dyn. Syst., and IMERL-Fac. Ing., CM-Fac. Cien.-Uruguay.

of open sets of diffeomorphisms whose elements display infinitely many sinks. In subsequent decades, other fundamental dynamic prototypes were found in this context, namely the so-called Hénon-like strange attractor ([BC], [MV]), and even infinitely many coexisting ones [C]. Even before these last results, Palis ([PT], [P1]) was already conjecturing that the presence of a homoclinic bifurcation (homoclinic tangency in the case of surfaces) is a very common phenomenon in the complement of the closure of the $\Omega$-stable ones (or hyperbolic ones). In fact, if the conjecture is proved to be true, then homoclinic bifurcation would certainly play a central role in the global understanding of the space of dynamics, for it would imply that each of these bifurcation phenomena is dense in the complement of the closure of the $\Omega$-stable ones. More precisely, the conjecture is:

CONJECTURE (Palis). *Every $f \in \mathrm{Diff}^r(M), r \geq 1$, can be $C^r$-approximated by a diffeomorphism exhibiting a homoclinic bifurcation or by one which is (essentially) hyperbolic.*

Here, hyperbolic means a diffeomorphism such that its limit set is hyperbolic, and essentially hyperbolic refers to a diffeomorphism having finitely many hyperbolic attractors whose basin of attraction covers a set of total probability (Lebesgue). When $M$ is a surface we can replace, in the above conjecture, the term "homoclinic bifurcation" by "homoclinic tangency." The aim of this work is to prove this conjecture in case $M$ is a surface and $r = 1$.

THEOREM A. *Let $M^2$ be a two dimensional compact manifold and let $f \in \mathrm{Diff}^1(M^2)$. Then, $f$ can be $C^1$-approximated by a diffeomorphism exhibiting a homoclinic tangency or by an Axiom A diffeomorphism.*

We recall that the stable and unstable sets
$$W^{\mathrm{s}}(p, f) = \{y \in M : \mathrm{dist}(f^n(y), f^n(p)) \to 0 \text{ as } n \to \infty\},$$
$$W^{\mathrm{u}}(p, f) = \{y \in M : \mathrm{dist}(f^n(y), f^n(p)) \to 0 \text{ as } n \to -\infty\}$$
are $C^r$-injectively immersed submanifolds when $p$ is a hyperbolic periodic point of $f$. A point of intersection of these manifolds is called a homoclinic point. We say that a diffeomorphism exhibits a homoclinic tangency if for some hyperbolic periodic point of it, the manifolds defined above are tangent at some point.

Also, a set $\Lambda$ is called hyperbolic for $f$ if it is compact, $f$-invariant and the tangent bundle $T_\Lambda M$ can be decomposed as $T_\Lambda M = E^{\mathrm{s}} \oplus E^{\mathrm{u}}$ invariant under $Df$ and there exist $C > 0$ and $0 < \lambda < 1$ such that
$$\|Df^n_{/E^{\mathrm{s}}(x)}\| \leq C\lambda^n$$
and
$$\|Df^{-n}_{/E^{\mathrm{u}}(x)}\| \leq C\lambda^n$$

for all $x \in \Lambda$ and for every positive integer $n$. Moreover, a diffeomorphism is called Axiom A, if the nonwandering set is hyperbolic and if it is the closure of the periodic points. For the nice dynamic properties of an Axiom A diffeomorphism we refer to [B] and [Sh].

Let us discuss some consequences of Theorem A. In his improved version of a result by Birkhoff, Smale [S2] showed that in the presence of a transversal homoclinic orbit (the stable and unstable manifolds meet transversally at some point) there exist very rich and highly developed forms of recurrence (i.e. nontrivial hyperbolic sets). Thus, a natural question is if such a phenomenon is common or abundant in the complement of the closure of the "simple" systems, i.e. Morse-Smale systems, defined precisely by having a "simple" nonwandering set, formed by only finitely many orbits, all of them periodic and hyperbolic and having all their stable and unstable manifolds transverse to each other. Moreover, it is well known that Morse-Smale diffeomorphisms always have zero entropy, but it is not known whether a diffeomorphism having zero entropy can be approximated by a Morse-Smale one. Indeed, it is an old question whether the set of diffeomorphisms having zero entropy is contained in the closure of the Morse-Smale ones. We can give a positive answer to the first question and a partial one to the last one for surface diffeomorphisms in the $C^1$ topology:

COROLLARY. *Let* M-S *be the set of Morse-Smale diffeomorphisms, and denote by* Cl(M-S) *its closure and consider* $\mathcal{U} = \mathrm{Diff}^1(\mathrm{M}^2) - \mathrm{Cl(M\text{-}S)}$. *Then, there exists an open and dense set* $\mathcal{R}$ *in* $\mathcal{U}$ *such that every* $f \in \mathcal{R}$ *has a transversal homoclinic orbit. In particular, the closure of the interior of the set formed by the diffeomorphisms having zero entropy, is equal to* Cl(M-S).

This corollary follows immediately from Theorem A. In fact, since the set of diffeomorphisms having a transversal homoclinic orbit is open and dense in the set of those having homoclinic orbits, we only have to prove that in the complement of the closure of the Morse-Smale dynamical systems, there exists a dense subset exhibiting a homoclinic orbit. In other words, if a diffeomorphism $f$ cannot be approximated by a Morse-Smale one, then it can be approximated by one having a homoclinic orbit. Thus, assume that $f$ is not approximated by a Morse-Smale diffeomorphism. By Theorem A, such a diffeomorphism can be approximated by one having a homoclinic tangency (and the result holds) or by an Axiom A one. Then, if $f$ is approximated by Axiom A diffeomorphisms, and since no one of these is approximated by Morse-Smale ones, they must have homoclinic orbits. The second part of the corollary follows from the fact that a diffeomorphism having a transversal homoclinic orbit always has nonzero entropy.

Concerning the proof of Theorem A, we will show that if a diffeomorphism $f$ cannot be approximated by a diffeomorphism having a homoclinic tangency, then it can be approximated by an Axiom A one. This will be done in two steps: first, using arguments developed by Mañé in the proof of the stability conjecture [M1], we will show that it is possible to find a continuous splitting with certain special properties on "almost" the whole nonwandering set of $g$, for $g$ in some dense open subset of a neighborhood of $f$. Secondly, we will find $g$ near $f$ having this continuous splitting hyperbolic.

Let us be more precise. An $f$-invariant set $\Lambda$ is said to have a dominated splitting if we can decompose its tangent bundle in two invariant subbundles $T_\Lambda M = E \oplus F$, such that there exist $C > 0$ and $0 < \lambda < 1$ with the following property:

$$\|Df^n_{/E(x)}\|\|Df^{-n}_{/F(f^n(x))}\| \leq C\lambda^n, \text{ for all } x \in \Lambda, n \geq 0.$$

If $f$ cannot be approximated by a diffeomorphism having a homoclinic tangency, we will prove that the angle between the stable and unstable subspaces of every hyperbolic periodic point for every $g$ near $f$ is bounded away from zero, and this implies the existence of a dominated splitting on the nonwandering set of $g$. Unfortunately, we cannot expect to proceed with the same arguments as in the proof of the stability conjecture, which consists in showing that, if the decomposition is not hyperbolic, we can perturb the map to find a nonhyperbolic periodic point. This is a contradiction because in a neighborhood of a stable diffeomorphism all periodic points are hyperbolic. In our case, dominated splitting can coexist with nonhyperbolic periodic points, and we do not know whether in a neighborhood of $f$ there is or is not a dense subset of diffeomorphisms having a nonhyperbolic periodic point. Thus, we shall pursue a different argument: we will find $g$ near $f$ with a certain property which implies that the dominated splitting is hyperbolic. Surprisingly, this "special" property turns out to be smoothness (together with other well-known generic properties). In fact, inspired by a result of Mañé in one dimensional dynamics [M2], we have the following theorem:

THEOREM B. *Let $f$ be a $C^2$-diffeomorphism on a compact surface, $\Lambda$ a compact $f$-invariant set having a dominated splitting $T_{/\Lambda}M = E \oplus F$. Assume that all the periodic points in $\Lambda$ are hyperbolic of saddle type. Then $\Lambda = \Lambda_1 \cup \Lambda_2$, where $\Lambda_1$ is hyperbolic and $\Lambda_2$ consists of a finite union of periodic simple closed curves $\mathcal{C}_1, \ldots \mathcal{C}_n$, normally hyperbolic and such that $f^{m_i} : \mathcal{C}_i \to \mathcal{C}_i$ is conjugated to an irrational rotation ($m_i$ denotes the period of $\mathcal{C}_i$).*

We observe that a weaker version of Theorem B appeared in [A]. However, not only certain parts of the proof are unclear to us, but also it is not useful in our context. For this reason, we had to improve considerably the arguments in [M2], in order to work in the two dimensional case. We also want to point out that Mañé orally announced in the late eighties that he and Araujo had proved Palis's conjecture in the $C^1$ topology. Nevertheless, this has never been shown and our paper aims to provide a proof of it.

The proof of Theorem A is presented in Section 2 if we assume Theorem B. We then show that Theorem B is true in Section 3.

*Acknowledgments.* This is the second author's thesis at IMPA under the guidance of Jacob Palis. We are very grateful to him for many valuable commentaries and for all his encouragement. Also, to Marcelo Viana for many useful conversations.

## 2. Proof of Theorem A

As mentioned in the introduction, let $M$ be a smooth compact surface endowed with some Riemannian metric and let $\text{Diff}^1(M)$ be the set of $C^1$ diffeomorphisms of $M$ endowed with the $C^1$ uniform topology. In this section, we shall prove Theorem A assuming that Theorem B is true. So, for any $f \in \text{Diff}^1(M)$ we shall show that either $f$ can be approximated by a diffeomorphism exhibiting a homoclinic tangency or by an Axiom A one.

Let $E, F$ be two subspaces of $\mathbb{R}^2$, $\mathbb{R}^2 = E \oplus F$. The angle$(E, F)$ is defined as the norm of the operator $L : E \to E^\perp$ such that

$$F = \text{graph}(L) = \{u + Lu : u \in E\}.$$

For $f \in \text{Diff}^1(M)$, denote by $\text{Per}_h(f)$ the set of hyperbolic periodic points of saddle type of $f$. Moreover, let $E_p^s(f), E_p^u(f)$ be the stable and unstable subspaces of a hyperbolic periodic point of $f$.

LEMMA 2.0.1. *Let $f \in \text{Diff}^1(M)$ and assume that there exist $\gamma > 0$ and a neighborhood $\mathcal{U}(f)$ of $f$ such that for every $g \in \mathcal{U}(f)$ and every $p \in \text{Per}_h(g)$,*

$$\text{angle}(E_p^s(g), E_p^u(g)) > \gamma$$

*holds. Then, the closure $\text{Cl}(\text{Per}_h(f))$ has a dominated splitting.*

The proof will be given in subsection 2.1.

For $f \in \text{Diff}^1(M)$, denote by $P_0(f)$ the set of sinks of $f$ and by $F_0(f)$ the set of sources. Let $\Omega(f)$ be the nonwandering set and define

$$\Omega_0(f) = \Omega(f) - (P_0(f) \cup F_0(f)).$$

Notice that $\Omega_0(f)$ is compact.

Let $\mathcal{U}$ be the complement of the closure of the diffeomorphisms exhibiting a homoclinic tangency, i.e.,

$$\mathcal{U} = \text{Diff}^1(M) - \text{Cl}(\{f \in \text{Diff}^1(M) : f \text{ exhibits a homoclinic tangency}\}).$$

Now, we have:

LEMMA 2.0.2. *There exists $\mathcal{U}_1 \subset \mathcal{U}$ which is open and dense in $\mathcal{U}$ such that, for every $g \in \mathcal{U}_1$, $\Omega_0(g)$ has a dominated splitting.*

The proof will be given in subsection 2.2.

Now, we can conclude the proof of Theorem A. Let $f \in \text{Diff}^1(M)$ be any diffeomorphism. We must show that it can be $C^1$-approximated by a diffeomorphism exhibiting a homoclinic tangency or by one which is Axiom A. Let us assume that $f$ cannot be $C^1$-approximated by a diffeomorphism exhibiting a homoclinic tangency. Then $f \in \mathcal{U}$. We can take a $C^2$ Kupka-Smale diffeomorphism $g$ (i.e., each periodic point is hyperbolic and all stable and unstable manifolds are in general position) arbitrarily close to $f$ and such that $g \in \mathcal{U}_1$ (such diffeomorphisms form a residual subset of $\text{Diff}^r(M)$, for all $r \geq 1$ [K], [S3]). Thus, by the previous lemma, $\Omega_0(g)$ has a dominated splitting. Moreover, the periodic points of $g$ in $\Omega_0(g)$ are hyperbolic of saddle type.

By Theorem B, $\Omega_0(g)$ is the union of a hyperbolic set and a finite union of normally hyperbolic periodic curves conjugated to an irrational rotation. Since the existence of such curves is not generic (even in the $C^2$ topology), and $\mathcal{U}_1$ is open, we can assume that $g$ does not have such invariant curves.

Thus, $\Omega_0(g)$ is hyperbolic. We want to show that $\Omega(g)$ is also hyperbolic. Since $\Omega(g) = \Omega_0(g) \cup P_0(g) \cup F_0(g)$, we only need to show that the set of sinks is a finite set and so is the set of sources.

Assume that $\#P_0(g) = \infty$. Notice that $\text{Cl}(P_0(g)) - P_0(g)$ is contained in $\Omega_0(g)$ which is hyperbolic. Moreover, there exists a neighborhood $U$ of $\Omega_0(g)$ such that the maximal invariant set in this neighborhood is hyperbolic. Since $\text{Cl}(P_0(g)) - P_0(g)$ is contained in $\Omega_0(g)$, we conclude that all the sinks (except, perhaps, a finite number of them) are contained in $U$, a contradiction to the hyperbolicity of the maximal invariant subset of $U$. Thus, $\#P_0(g)$ is finite. The same arguments apply to $F_0(g)$.

We conclude then that $\Omega(g)$ is hyperbolic. As $g$ is Kupka-Smale and $M$ has dimension two, it is not difficult to see that $g$ is an Axiom A diffeomorphism. In fact, since $\Omega(g)$ is hyperbolic, the limit set $L(f)$ is hyperbolic. By a result of Newhouse [N4], the periodic points are dense in $L(f)$ and we can do a spectral decomposition of $L(f)$. To prove that $\Omega(f) = L(f)$ we only have to check the no-cycle condition on $L(f)$. Since $M$ has dimension two, a cycle can only be formed by basic sets of saddle type (or index one). From the fact that $g$ is a Kupka-Smale diffeomorphism we can conclude that all the intersections of this cycle are transversal. Finally, the last two assertions imply that the cycle does not exist. This completes the proof of Theorem A.

2.1. *Proof of Lemma* 2.0.1. First we need the next lemma due to Franks. The version we shall use here is slightly more general but the proof in [F] can be easily extended to include this statement.

LEMMA 2.1.1. *Let $f \in \mathrm{Diff}^1(M)$ and let $\mathcal{U}(f)$ be any neighborhood of $f$. Then, there exist $\varepsilon > 0$ and $\mathcal{U}_0(f) \subset \mathcal{U}(f)$ such that given $g \in \mathcal{U}_0(f)$, a finite set $\{x_1, \ldots, x_N\}$, a neighborhood $U$ of $\{x_1, \ldots, x_N\}$ and linear maps $L_i : T_{x_i}M \to T_{g(x_i)}M$ such that $\|L_i - D_{x_i}g\| \leq \varepsilon$ for all $1 \leq i \leq N$, then there exists $\tilde{g} \in \mathcal{U}(f)$ such that $\tilde{g}(x) = g(x)$ if $x \in \{x_1, \ldots, x_N\} \cup M - U$ and $D_{x_i}\tilde{g} = L_i$ for all $1 \leq i \leq N$.*

Now, let $\mathcal{U}(f)$ be a neighborhood as in the statement of Lemma 2.0.1 and let $\varepsilon > 0$ and $\mathcal{U}_0(f)$ be as in the previous lemma.

We claim that there exists $\delta > 0$ such that for every $p \in \mathrm{Per}_h(g), g \in \mathcal{U}_0(f)$ we have either

$$|\lambda_p| < (1 - \delta)^n$$

or

$$|\sigma_p| > (1 + \delta)^n$$

where $\lambda_p$ and $\sigma_p$ are the eigenvalues of $Dg^n : T_pM \to T_pM$ ($n$ is the period of $p$), $|\lambda_p| < 1 < |\sigma_p|$.

Let $C = \sup\{\|Dg\| : g \in \mathcal{U}_0(f)\}$ and consider $\varepsilon' < \frac{\varepsilon}{C}$.

To prove this claim, let us observe first that there exists $\beta_0 = \beta_0(\varepsilon')$, $0 < \beta_0 < 1$, such that for every $\beta', \beta'' < \beta_0$ and for every vector $u, v \in T_xM$ such that $\mathrm{angle}(u, v) > \gamma$ then, if we take $T : T_xM \to T_xM$ satisfying

$$T_{/<u>} = (1 - \beta')\mathrm{Id}, \quad T_{/<v>} = (1 + \beta'')\mathrm{Id}$$

we have that

$$\|T - \mathrm{Id}\| < \frac{\varepsilon'}{2}.$$

Assume that the claim is not true. Then, there exists a sequence $\delta_m \to 0$, $g_m \in \mathcal{U}_0$ and $p_m \in \mathrm{Per}_h(g_m)$ of period $n_m$ such that

$$(1 - \delta_m)^{n_m} < |\lambda_{p_m}| < 1$$

and

$$1 < |\sigma_{p_m}| < (1 + \delta_m)^{n_m}.$$

For simplicity, assume that $\lambda_{p_m}$ and $\sigma_{p_m}$ are positive. Take $m$ big enough such that $\delta_m < \beta_0$, and set $p = p_m, n = n_m, g = g_m$ to simplify the notation.

Let $\beta' = 1 - \lambda_p^{\frac{1}{n}}$ and $\beta'' = \sigma_p^{\frac{1}{n}} - 1$. For every $0 \leq i \leq n - 1$ define $T_i : T_{g^i(p)}M \to T_{g^i(p)}M$ in such a way that

$$T_{i/E^s_{g^i(p)}} = (1 - \beta')\mathrm{Id}$$

and

$$T_{i/E^u_{g^i(p)}} = (1 + \beta'')\mathrm{Id}.$$

Moreover, take two vectors $u_0, v_0 \in T_pM$, $\text{angle}(u_0, v_0) < \gamma$, and $\beta_1 > 0$ and define $S: T_pM \to T_pM$ such that

$$S_{/\langle u_0 \rangle} = (1 - \beta_1)\text{Id}, \ S_{/\langle v_0 \rangle} = (1 + \beta_1)\text{Id}.$$

Observe that, if $\beta_1$ is small enough, $\|S - \text{Id}\|\|T_0\| < \frac{\varepsilon'}{2}$.

Now, for $0 \leq i \leq n - 2$ define

$$L_i: T_{g^i(p)}M \to T_{g^{i+1}(p)}M \text{ by } L_i = T_{i+1} \circ Dg_{g^i(p)}$$

and

$$L_{n-1} = S \circ T_0 \circ Dg_{g^{n-1}(p)}.$$

Notice that

$$\|L_i - Dg_{g^i(p)}\| < \varepsilon \text{ for all } i = 0, 1, \ldots n - 1.$$

In fact, for $0 \leq i \leq n - 2$,

$$\begin{aligned}
\|L_i - Dg_{g^i(p)}\| &= \|T_{i+1} \circ Dg_{g^i(p)} - Dg_{g^i(p)}\| \\
&\leq \|T_{i+1} - \text{Id}\|\|Dg_{g^i(p)}\| < \frac{\varepsilon'}{2}C < \varepsilon
\end{aligned}$$

and

$$\begin{aligned}
\|L_{n-1} - Dg_{g^{n-1}(p)}\| &= \|S \circ T_0 \circ Dg_{g^{n-1}(p)} - Dg_{g^{n-1}(p)}\| \\
&\leq \|S \circ T_0 - \text{Id}\|\|Dg_{g^{n-1}(p)}\| \\
&\leq \left(\|S - \text{Id}\|\|T_0\| + \|T_0 - \text{Id}\|\right)C < \frac{\varepsilon}{2} + \frac{\varepsilon}{2} = \varepsilon.
\end{aligned}$$

Then, by Lemma 2.1.1, there exists $\tilde{g} \in \mathcal{U}(f)$ such that $p \in \text{Per}(\tilde{g})$ and $D\tilde{g}_{\tilde{g}^i(p)} = L_i$. Since $D\tilde{g}_p^n = L_{n-1} \circ \cdots \circ L_0 = S$, $p \in \text{Per}_h(\tilde{g})$. Moreover $E_p^s(\tilde{g}) = \langle u_0 \rangle$, $E_p^u(\tilde{g}) = \langle v_0 \rangle$ and so

$$\text{angle}(E_p^s(\tilde{g}), E_p^u(\tilde{g})) < \gamma,$$

a contradiction. This proves our claim.

For the rest of the proof we follow [M3]. To show the existence of a dominated splitting on $\text{Cl}(\text{Per}_h(f))$, it is enough to prove that the decomposition in stable and unstable spaces of the periodic points on $\text{Per}_h(f)$ is a dominated splitting, due to the fact that this dominated splitting can be easily extended to the closure having the same property.

Let us prove then, the dominated splitting on $\text{Per}_h(f)$. Observe that it is enough to show the existence of a positive integer $m_1$ such that for every $p \in \text{Per}_h(f)$ there exists $1 \leq m \leq m_1$ such that

$$\|Df_{/E_p^s}^m\|\|Df_{/E_{f^m(p)}^u}^{-m}\| < \frac{1}{2}.$$

Arguing by contradiction, assume that such $m_1$ does not exist. Then, we can find a sequence $p_n \in \text{Per}_h(f)$ satisfying

$$\|Df^m_{/E^s_{p_n}}\|\|Df^{-m}_{/E^u_{f^m(p_n)}}\| \geq \frac{1}{2}$$

for every $0 < m \leq n$. Due to our claim the periods of the periodic points $p_n$ must be unbounded. Otherwise, if the periods of $p_n$ are bounded, we can assume (taking a subsequence if necessary) that all have the same period $m_0$, and $p_n$ converges to a periodic point $p$ of period at most $m_0$. Observe that the point could not be hyperbolic, but has two real eigenvalues, one of them having modulus different from 1. Moreover, we can assume that $E^s_{p_n}$ converges to an invariant space $E_p$ and $E^u_{p_n}$ converges to an invariant space $F$ (which is necessarily different from $E$). Then, for all $m \geq 0$,

$$\|Df^m_{/E_p}\|\|Df^{-m}_{/F_{f^m(p)}}\| \geq \frac{1}{2};$$

but this is impossible since one eigenvalue is different from 1. Thus, the periods of $p_n$ are unbounded.

We can assume without loss of generality that for every $p_n$ we have $\lambda_n < (1-\delta)^{m_n}$, where $m_n$ is the period of $p_n$. Now, take $\varepsilon_0 > 0$ satisfying $(2\varepsilon_0 + \varepsilon_0^2)C \leq \varepsilon$, $\varepsilon_1 > 0$, and $m$ such that

$$(1 + \varepsilon_1)(1-\delta) < 1,$$
$$\varepsilon_1 \leq \frac{\gamma}{1+\gamma}\varepsilon_0$$

and

$$(1 + \varepsilon_1)^m \geq 4 + \frac{2}{\gamma}.$$

Since the periods of the periodic points are unbounded, we can choose $p_n$ such that $m_n > m$. Set $p = p_n$ and $n_0 = m_n$. Take $v \in E^u_p$ and $w \in E^s_p$ with $\|v\| = \|w\| = 1$. Hence

$$\frac{1}{2}\|Df^m v\| \leq \|Df^m w\|.$$

Take a linear map $L : E^u_p \to E^s_p$ satisfying

$$Lv = \varepsilon_1 w \text{ and } \|L\| = \varepsilon_1.$$

Define $\tilde{L} : E^u_p \to E^s_p$ by

$$\tilde{L} = (1 + \varepsilon_1)^{n_0}(Df^{n_0}_{/E^s_p}) \circ L \circ Df^{-n_0}_{/E^u_p}.$$

Observe that

$$\|\tilde{L}\| \leq (1+\varepsilon_1)^{n_0}(1-\delta)^{n_0}\|L\| \leq \varepsilon_1.$$

Define

$$G = \{u + Lu : u \in E^u_p\}$$

and
$$\tilde{G} = \{u + \tilde{L}u : u \in E^{\mathrm{u}}_p\}$$

and take linear maps $P, S$ from $T_pM$ to itself such that
$$\begin{aligned} P_{/E^{\mathrm{s}}_p} &= 0, \\ (\mathrm{Id} + P).E^{\mathrm{u}}_p &= G, \\ S_{/E^{\mathrm{s}}_p} &= 0 \end{aligned}$$

and
$$(\mathrm{Id} + S).\tilde{G} = E^{\mathrm{u}}_p.$$

It is possible to choose these maps satisfying (see Lemma II.10 of [M3])
$$\|P\| \leq \frac{1+\gamma}{\gamma} \|L\| \leq \frac{1+\gamma}{\gamma}\varepsilon_1 \leq \varepsilon_0$$

and
$$\|S\| \leq \frac{1+\gamma}{\gamma} \|\tilde{L}\| \leq \frac{1+\gamma}{\gamma}\varepsilon_1 \leq \varepsilon_0.$$

Finally take linear maps $T_j : T_{f^j(p)}M \to T_{f^j(p)}M$, $0 \leq j \leq n_0 - 1$, satisfying
$$T_{j/E^{\mathrm{s}}_{f^j(p)}} = \varepsilon_1 \mathrm{Id}$$

and
$$T_{j/E^{\mathrm{u}}_{f^j(p)}} = 0.$$

The norm of $T_j$ can be estimated by
$$\|T_j\| \leq \frac{1+\gamma}{\gamma}\varepsilon_1 \leq \varepsilon_0.$$

With these maps we shall construct a perturbation of $Df$ along the orbit of $p$. Let $L_0 = (\mathrm{Id} + T_1) \circ Df \circ (\mathrm{Id} + P)$ and for $1 \leq j \leq n_0 - 2$, set
$$L_j = (\mathrm{Id} + T_{j+1}) \circ Df$$

and
$$L_{n_0-1} = (\mathrm{Id} + S) \circ (\mathrm{Id} + T_0) \circ Df.$$

It follows that
$$\|L_j - Df_{f^j(p)}\| \leq \varepsilon, \text{ for } 0 \leq j \leq n_0 - 1.$$

Then, Lemma 2.1.1 implies that there exists $g \in \mathcal{U}(f)$ such that $p$ is a periodic point of $g$ and $Dg_{g^j(p)} = L_j$. Therefore $\beta = \mathrm{angle}(E^{\mathrm{s}}_{g^m(p)}(g), E^{\mathrm{u}}_{g^m(p)}(g)) > \gamma$. We shall show that this is a contradiction. For simplicity, denote $E^{\mathrm{s}}_j(f) = E^{\mathrm{s}}_{f^j(p)}(f)$, $E^{\mathrm{u}}_j(f) = E^{\mathrm{u}}_{f^j(p)}(f)$.

Observe that $E^{\mathrm{s}}_j(g) = E^{\mathrm{s}}_j(f)$ because
$$Dg_{/E^{\mathrm{s}}_j(f)} = L_j = (1+\varepsilon_1)Df_{/E^{\mathrm{s}}_j(f)}$$

and
$$(1+\varepsilon_1)(1-\delta) < 1.$$

On the other hand we have
$$Dg^{n_0}.E_0^{\mathrm{u}}(f) = L_{n_0-1} \circ \cdots \circ L_0.E_0^{\mathrm{u}}(f) = E_0^{\mathrm{u}}(f)$$
and
$$\|Dg^{n_0}_{/E_0^{\mathrm{u}}(f)}\| \geq \|Df^{n_0}_{/E_0^{\mathrm{u}}(f)}\| > 1.$$

Hence, $p \in \mathrm{Per}_h(g)$ and
$$E_0^{\mathrm{u}}(g) = E_0^{\mathrm{u}}(f).$$

Now, let us estimate the angle $\beta$ between $E_m^{\mathrm{s}}(g)$ and $E_m^{\mathrm{u}}(g)$. Since $v \in E_0^{\mathrm{u}}(g)$ and $w \in E_0^{\mathrm{s}}(g)$,
$$u_1 = Dg^m.v \in E_m^{\mathrm{u}}(g)$$
and
$$u_2 = Dg^m.w \in E_m^{\mathrm{s}}(g).$$

Moreover
$$u_1 = Df^m.v + (1+\varepsilon_1)^m Df^m.w$$
and
$$u_2 = (1+\varepsilon_1)^m Df^m.w.$$

Hence (see Lemma II.10 of [M3]),
$$\|Df^m.v\| = \|u_1 - u_2\| \geq \frac{\beta}{1+\beta}\|u_1\|$$
$$\geq \frac{\beta}{1+\beta}\Big|(1+\varepsilon_1)^m\|Df^m.w\| - \|Df^m.v\|\Big|.$$

Therefore
$$\frac{1+\beta}{\beta} \geq \frac{(1+\varepsilon_1)^m}{2} - 1.$$

Thus
$$\beta \leq \frac{2}{(1+\varepsilon_1)^m - 4} \leq \gamma,$$
a contradiction, completing the proof of the lemma.

2.2. *Proof of Lemma* 2.0.2. The following lemma represents the key for the relationship between homoclinic tangencies and the angle of the stable and unstable subspaces of periodic points.

LEMMA 2.2.1. *Let $f$ be any diffeomorphism and fix $\varepsilon > 0$. Assume that there exists a hyperbolic periodic point $p$ for $f$ with eigenvalues $\lambda$, $\sigma$, $|\lambda| < 1 < |\sigma|$, $|\lambda\sigma| < 1$ so that $\gamma = \mathrm{angle}(E_p^{\mathrm{s}}, E_p^{\mathrm{u}})$ satisfies*
$$\gamma < \frac{|\sigma|-1}{|\sigma|+1}\frac{\varepsilon}{2}.$$

*Then, there exist a diffeomorphism $g, \varepsilon$-$C^1$-close to $f$ having $p$ as a hyperbolic periodic point, and a homoclinic tangency associated to $p$.*

*Proof.* Take a very small neighborhood $B(p)$ of $p$ satisfying:

1. $\mathcal{O}(p) \cap B(p) = \{p\}$ ($\mathcal{O}(p)$ is the orbit of $p$).

2. $f^j(B(p)) \cap B(p) = \emptyset$ for $1 \leq |j| < n$ where $n$ is the period of $p$.

3. $\exp_p : B(0,r) \to B(p)$ is a diffeomorphism where $B(0,r) \subset T_pM$ is a ball of radius $r$ at the origin.

We shall make a perturbation on $T_pM$ and then we compose with the exponential map to pass to the manifold. We shall assume that $0 < \lambda < 1 < \sigma$, the other cases are similar.

Let us put coordinates on $T_pM = E_p^{\mathrm{u}} \oplus E_p^{\mathrm{u}\perp}$. Assume, without loss of generality, that $E_p^{\mathrm{s}} = \{v + \gamma v : v \in E_p^{\mathrm{u}}\}$.

Let $\phi : \mathbb{R} \to \mathbb{R}$ and $\psi : \mathbb{R} \to \mathbb{R}$ be two $C^\infty$ maps such that

- $0 \leq \phi \leq \varepsilon, |\phi'| \leq \varepsilon, \phi(0) = \varepsilon$ and $\phi(x) = 0$ for $|x| > 2$.

- $0 \leq \psi \leq 1, |\psi'| \leq 1, \psi(0) = 1$ and $\psi(x) = 0$ for $|x| > 2$.

Now, for every $a$ consider the map $\Phi_a : T_pM \to T_pM$ defined by

$$\Phi_a(x,y) = (x,y) + (0, a\phi(\frac{x}{a})\psi(\frac{y}{a})).$$

It is not difficult to see that

1. $\Phi_a$ is a diffeomorphism;

2. $\Phi_a$ is $\varepsilon$-$C^1$-close to the identity;

3. $\Phi_a(0,0) = (0, a\varepsilon)$;

4. $\Phi_a = \mathrm{id}$ if $|x| > 2a, |y| > 2a$.

For simplicity, we shall assume that the stable and unstable subspaces of $Df_p^n : T_pM \to T_pM$ are mapped onto the local stable and unstable manifolds of $p$ under $\exp_p$ restricted to $B(0,r)$. This may be assumed because the stable and unstable manifolds (under $\exp_p^{-1}$) are $C^1$ close to the stable and unstable subspaces and tangent at the origin, and we shall do the perturbation arbitrarily close to the origin.

Let $x_1$ be any point on $E_p^{\mathrm{u}}$ arbitrarily close to the origin, and consider $(x_1, \sigma x_1)$ a fundamental domain in $E_p^{\mathrm{u}}$. Take $x_0 = \frac{\sigma x_1 + x_1}{2}$ as the middle point in this fundamental domain, and let $a = \frac{(\sigma-1)x_1}{4}$. For $0 \leq t \leq 1$ consider the family $\tilde{\Phi}_{ta} : T_pM \to T_pM$ defined by

$$\tilde{\Phi}_{ta}(x,y) = \Phi_{ta}(x - x_0, y).$$

Observe that $\tilde{\Phi}_{ta}$ is equal to the identity if $|x - x_0| > t\frac{(\sigma-1)x_1}{2}$, $|y| > t\frac{(\sigma-1)x_1}{2}$.

Now, for $t = 1$, we have $\tilde{\Phi}_a(x_0, 0) = (0, a\varepsilon)$ and

$$a\varepsilon = \frac{(\sigma-1)x_1}{4}\varepsilon = \frac{\sigma-1}{\sigma+1}x_1\frac{\varepsilon}{2}\frac{\sigma+1}{2} > \gamma x_1\frac{\sigma+1}{2} = \gamma x_0.$$

Since $\text{angle}(E_p^s, E_p^u) = \gamma$, that is, $E_p^s$ is the graph of the map $L : E_p^u \to E_p^{u\perp}$, $Lv = \gamma v$, we conclude that the curve

$$\{\tilde{\Phi}_a(x, 0), \ x_1 \leq x \leq \sigma x_1\}$$

intersects $E_p^s$. Hence, for some $t_0 \leq 1$ we have that the curve

$$\{\tilde{\Phi}_{t_0 a}(x, 0), \ x_1 \leq x \leq \sigma x_1\}$$

is tangent to $E_p^s$.

Let $D = \{(x, y) \in T_p M : |x - x_0| \leq \frac{\sigma-1}{2}x_1, |y| \leq \frac{\sigma-1}{2}x_1\}$. Notice that $\tilde{\Phi}_{ta} = \text{id}$ outside $D$ and moreover $E_p^s \cap D$ is contained in a fundamental domain of $E_p^s$ (remember that $\lambda\sigma < 1$).

Consider the map $g : M \to M$, $g = f$, on $M - f^{-1}(B(p))$ and if $x \in f^{-1}(B(p))$ define

$$g(x) = \exp_p \circ \tilde{\Phi}_{t_0 a} \circ \exp_p^{-1}(f(x)).$$

It is easy to see that $g$ is a diffeomorphism $\varepsilon$-$C^1$-close to $f$, and $p$ is a hyperbolic periodic point of $g$ (observe that $0 \notin D$).

To finish the proof of our lemma we only need to verify that the curve defined by $\exp_p\{(\tilde{\Phi}_{t_0 a}(x, 0)), x_1 \leq x \leq \sigma x_1\}$ is in $W^u(p, g)$ and the curve $\exp_p(E_p^s \cap D)$ is in $W^s(p, g)$, since we already know that they are tangent.

- Let $y \in \exp_p(E_p^s \cap D)$. Since $f^n(y) \cap \exp(D) = \emptyset$ for every positive $n$, we have that $g^n(y) = f^n(y)$, and since $y \in W^s(p, f)$ we conclude also that $y \in W^s(p, g)$.

- Let $y \in \exp_p\{(\tilde{\Phi}_{t_0 a}(x, 0)), x_1 \leq x \leq \sigma x_1\}$. Then the point

$$y_1 = (\exp_p \circ \tilde{\Phi}_{t_0 a} \circ \exp_p^{-1})^{-1}(y)$$

is in $W^u(p, f)$. Moreover $f^{-n}(y_1) \cap \exp(D) = \emptyset$ for all positive $n$. Thus, $g^{-n}(y) = f^{-n}(y_1)$, implying $y \in W^u(p, g)$.

This completes the proof of the lemma. □

Remember that $\mathcal{U}$ is the complement of the closure of the diffeomorphisms exhibiting a homoclinic tangency.

LEMMA 2.2.2. *Let $f \in \mathcal{U}$ and assume that $f$ is Kupka-Smale. Then there exist a neighborhood $\mathcal{U}(f)$ of $f$ and $\gamma > 0$ such that for every $g \in \mathcal{U}(f)$ and every $p \in \text{Per}_h(g)$*

$$\text{angle}(E_p^s, E_p^u) > \gamma.$$

*Proof.* Let $\mathcal{U}_1(f) \subset \mathcal{U}_0(f) \subset \mathcal{U}$ be two neighborhoods of $f$ and take $\varepsilon_1 > 0$ such that if $\tilde{g}$ is $\varepsilon_1$-$C^1$-close to $g \in \mathcal{U}_1(f)$ then $\tilde{g} \in \mathcal{U}_0(f)$. Moreover, consider $\mathcal{U}_2(f)$ and $\varepsilon > 0$ from Lemma 2.1.1 corresponding to $\mathcal{U}_1(f)$. Let $C = \sup\{\|Dg\| : g \in \mathcal{U}_0(f)\}$ and $\varepsilon' < \frac{\varepsilon}{C}$.

Now, arguing by contradiction, assume that there exist sequences $\gamma_n \to 0$, $g_n \to f$ and $p_n \in \text{Per}_h(g_n)$ such that

$$\text{angle}(E^s_{p_n}, E^u_{p_n}) = \gamma_n.$$

We can assume, without loss of generality, that $g_n \in \mathcal{U}_2(f)$, $\lambda_{p_n}\sigma_{p_n} < 1$, and

$$\gamma_n = \text{angle}(E^s_{p_n}, E^u_{p_n}) = \min\{\text{angle}(E^s_q, E^u_q) : q \text{ in the orbit of } p_n\}.$$

We shall show the existence of $\tilde{g} \in \mathcal{U}_1(f)$ and $p \in \text{Per}_h(\tilde{g})$ so that $\lambda_p \sigma_p < 1$ and such that

$$\text{angle}(E^s_p(\tilde{g}), E^u_p(\tilde{g})) = \gamma < \frac{(\sigma_p - 1)}{(\sigma_p + 1)} \frac{\varepsilon_1}{2}.$$

This leads to a contradiction because, by the previous lemma, we can find a diffeomorphism in $\mathcal{U}_0(f)$ exhibiting a homoclinic tangency.

So, if for some $n$,

$$\gamma_n < \frac{(\sigma_n - 1)}{(\sigma_n + 1)} \frac{\varepsilon_1}{2},$$

we can conclude the proof ($\sigma_n = \sigma_{p_n}, \lambda_n = \lambda_{p_n}$).

Hence, let us assume that this does not hold. Let $m_n$ be the period of $p_n$. Since $f$ is Kupka-Smale we have $m_n \to \infty$. Consider $\delta_n = \frac{\gamma_n \varepsilon'}{2}$ and for $0 \leq i \leq m_n - 1$ consider $T_i : T_{g_n^i(p_n)} \to T_{g_n^i(p_n)}$ such that

$$T_{i/E^s_{g_n^i(p_n)}} = (1 - \delta_n)\text{Id}$$

and

$$T_{i/E^u_{g_n^i(p_n)}} = (1 + \delta_n)\text{Id}.$$

It is not difficult to see that $\|T_i - \text{Id}\| < \varepsilon'$. Now, for $0 \leq i \leq m_n - 2$ consider

$$L_i : T_{g_n^i(p_n)} \to T_{g_n^{i+1}(p_n)}, \quad L_i = T_{i+1} \circ Dg_n(g_n^i(p_n))$$

and $L_{m_n-1} = T_0 \circ Dg_n(g_n^{m_n-1}(p_n))$. Now, we have $\|L_i - Dg_n(g_n^i(p_n))\| < \varepsilon$. Applying Lemma 2.1.1 we conclude the existence of $\tilde{g}_n \in \mathcal{U}_1(f)$ such that $p_n \in Per(\tilde{g}_n)$ and $D\tilde{g}_n(\tilde{g}_n^i(p_n)) = L_i$. Hence we have

$$\begin{aligned} p_n &\in \text{Per}_h(\tilde{g}_n), \\ \tilde{\gamma}_n &= \gamma_n, \\ \tilde{\lambda}_n &= (1 - \delta_n)^{m_n} \lambda_n, \end{aligned}$$

and

$$\tilde{\sigma}_n = (1 + \delta_n)^{m_n} \sigma_n.$$

Thus $\tilde{\lambda}_n \tilde{\sigma}_n < 1$ for all $n$. We shall prove that for $n$ arbitrarily large,

$$\frac{(\tilde{\sigma}_n - 1)}{(\tilde{\sigma}_n + 1)} \frac{\varepsilon_1}{2} > \tilde{\gamma}_n.$$

If the sequence $\tilde{\sigma}_n$ is unbounded, it follows immediately that we have the mentioned property (remember $\tilde{\gamma}_n = \gamma_n \to 0$). Therefore, assume for some $K$,

$$\tilde{\sigma}_n = (1 + \delta_n)^{m_n} \sigma_n \leq K \text{ for all } n.$$

Now

$$\begin{aligned}
\frac{(\tilde{\sigma}_n - 1)}{(\tilde{\sigma}_n + 1)} \frac{\varepsilon_1}{2} &= \frac{((1+\delta_n)^{m_n}\sigma_n - 1)}{((1+\delta_n)^{m_n}\sigma_n + 1)} \frac{\varepsilon_1}{2} \geq \frac{((1+m_n\delta_n)\sigma_n - 1)}{((1+\delta_n)^{m_n}\sigma_n + 1)} \frac{\varepsilon_1}{2} \\
&= \frac{(\sigma_n - 1)}{((1+\delta_n)^{m_n}\sigma_n + 1)} \frac{\varepsilon_1}{2} + \frac{(m_n\delta_n\sigma_n)}{((1+\delta_n)^{m_n}\sigma_n + 1)} \frac{\varepsilon_1}{2} \\
&\geq \frac{(m_n\delta_n\sigma_n)}{(K+1)} \frac{\varepsilon_1}{2} \geq \frac{(m_n\delta_n)}{(K+1)} \frac{\varepsilon_1}{2} = m_n \frac{\varepsilon'\varepsilon_1}{4(K+1)} \gamma_n > \gamma_n
\end{aligned}$$

where the last inequality is true if $n$ is big enough. $\square$

Let us now prove Lemma 2.0.2. Recall that a map $\Gamma$ which for each diffeomorphism $f : M \to M$ associates a compact subset $\Gamma(f)$ of $M$ is said to be lower semicontinuous if for every open set $U$ such that $U \cap \Gamma(f) \neq \emptyset$, there exists a neighborhood $\mathcal{U}(f)$ such that for every $g \in \mathcal{U}(f)$ we have $U \cap \Gamma(g) \neq \emptyset$; it is said to be upper semicontinuous if for every compact set $K$ such that $K \cap \Gamma(f) = \emptyset$, there exists a neighborhood $\mathcal{U}(f)$ with the property that for every $g \in \mathcal{U}(f)$ we have $K \cap \Gamma(g) = \emptyset$. And it is said to be continuous if it is lower and upper semicontinuous. It is well-known that a lower semicontinuous map $\Gamma$ has a residual set of points of continuity.

Let us take $\Gamma$ as $\Gamma(g) = \text{Cl}(\text{Per}_h(g))$. This map is lower semicontinuous since hyperbolic periodic points cannot be destroyed by any arbitrarily small perturbation. Thus, there exists a residual set $\mathcal{R}_1$ in $\text{Diff}^1(M)$ of points of continuity of the map $\Gamma$. Moreover, the set of Kupka-Smale diffeomorphisms forms a residual set $\mathcal{R}_2$, and there exists another residual $\mathcal{R}_3$ such that for every $g \in \mathcal{R}_3$ we have $\Omega(g) = \text{Cl}(\text{Per}(g))$ (see [Pu1] and [Pu2]). Let $\mathcal{R} = \mathcal{R}_1 \cap \mathcal{R}_2 \cap \mathcal{R}_3$ which is also a residual set, and let $\mathcal{R}_0 = \mathcal{R} \cap \mathcal{U}$ be a residual subset of $\mathcal{U}$ where

$$\mathcal{U} = \text{Diff}^1(M) - \text{Cl}(\{f \in \text{Diff}^1(M) : f \text{ exhibits a homoclinic tangency}\}).$$

Recall also that for any diffeomorphism $g$ we have defined

$$\Omega_0(g) = \Omega(g) - (P_0(g) \cup F_0(g))$$

where $P_0(g)$ is the set of sinks of $g$, and $F_0(g)$ the set of sources.

Now, taking $f \in \mathcal{R}_0$, we have

1. $\Omega(f) = \mathrm{Cl}(\mathrm{Per}(f))$,

2. $\mathrm{Per}(f) = P_0(f) \cup \mathrm{Per}_h(f) \cup F_0(f)$,

3. $f$ is a continuity point of the map $g \to \mathrm{Cl}(\mathrm{Per}_h(g))$.

We claim that $\Omega_0(f) = \mathrm{Cl}(\mathrm{Per}_h(f))$, for which we need the following result due to Pliss:

THEOREM 2.1 ([Pl]). *Let $f$ be a $C^1$ diffeomorphism and assume that $\#P_0(f) = \infty$. Then, given $\varepsilon$ there exist $g\varepsilon$-$C^1$-close to $f$ and $p \in P_0(f)$ such that $p$ is a nonhyperbolic periodic point of $g$.*

This is not exactly the statement contained in [Pl], but the proof can be easily extended to obtain this result. An analogous result for sources holds also.

We remark here that if we have a nonhyperbolic periodic point $p$ of a diffeomorphism $g$, then it is possible to find another diffeomorphism $g_1$ arbitrarily close to $g$ and a point $q$ arbitrarily close to $p$ such that $q$ is a hyperbolic periodic point of saddle type of $g_1$.

To prove our claim, notice that $\mathrm{Cl}(\mathrm{Per}_h(f)) \subset \Omega_0(f)$. Thus, we only have to show that $\Omega_0(f) \subset \mathrm{Cl}(\mathrm{Per}_h(f))$. Arguing by contradiction, assume that there exists a point $x \in \Omega_0(f)$, and $x \notin \mathrm{Cl}(\mathrm{Per}_h(f))$. Take a neighborhood $U$ of $\mathrm{Cl}(\mathrm{Per}_h(f))$, such that $x \notin \mathrm{Cl}(U)$. Since $f \in \mathcal{R}_0$, we conclude that there exists a neighborhood $\mathcal{U}(f)$ such that for every $g \in \mathcal{U}(f)$ we have $\mathrm{Cl}(\mathrm{Per}_h(g)) \subset U$. Now, the point $x$ must be accumulated by periodic points of $f$, and since $x \in \Omega_0(f)$ is not a periodic point, we conclude that we have a sequence of sinks or sources accumulating on $x$. We can apply the preceding theorem and remark, and we conclude that there exists $g$ arbitrarily close to $f$, say in $\mathcal{U}(f)$, and $q \in \mathrm{Per}_h(g)$ arbitrarily close to $x$, and so $\mathrm{Per}_h(g)$ is not contained in $U$, a contradiction. This proves the claim.

Thus, Lemma 2.0.1 and Lemma 2.2.2 imply that $\Omega_0(f)$ has dominated splitting $E \oplus F$. Take an admissible neighborhood $U$ of $\Omega_0(f)$, that is, there exists a neighborhood $\mathcal{U}(f)$ such that for every $g \in \mathcal{U}(f)$, the maximal invariant set of $g$ in $U$ also has a dominated splitting. We can assume also $\mathcal{U}(f) \subset \mathcal{U}$. Moreover, taking another neighborhood $U_0$ of $\Omega_0(f)$ such that $U_0 \subset \mathrm{Cl}(U_0) \subset U$, we can assume that for every $g \in \mathcal{U}(f)$ we have $\mathrm{Cl}(\mathrm{Per}_h(g)) \cap (M - U_0) = \emptyset$, since $f$ is a continuity point of the map $g \to \mathrm{Cl}(\mathrm{Per}_h(g))$.

Let us prove that for $g \in \mathcal{U}(f)$ we have $\Omega_0(g) \subset U$. Assume this is not true. Then, there exists $g \in \mathcal{U}(f)$ such that $\Omega_0(g) \cap (M - U) \neq \emptyset$. Take a point $x \in \Omega_0(g) \cap (M - U)$. The point $x$ is not periodic and it is in the nonwandering set $\Omega(g)$. We assert that the point $x \in \Omega_0(g) \cap (M - U)$ cannot be accumulated by sinks or sources of $g$. Indeed, if there exists a sequence $p_n$ of sinks (or sources) $p_n \to x$, applying again the preceding theorem and remark, we find a $\tilde{g}$ arbitrarily close to $g$ (say in $\mathcal{U}(f)$), having a hyperbolic periodic point of saddle type arbitrarily close to $x$, contradicting the continuity of the map $\Gamma$ at $f$. Now, by the $C^1$-closing lemma [Pu1] we can find $\tilde{g}$ close to $g$ (say in $\mathcal{U}(f)$) having $x$ as a periodic point. Moreover (perturbing a little if necessary) we can assume that $x$ is a hyperbolic periodic point of saddle type of $\tilde{g}$. Again, this is a contradiction with the continuity of the map $\Gamma$ at $f$. Thus, we have proved that for every $g \in \mathcal{U}(f), \Omega_0(g) \subset U$ and therefore $\Omega_0(g)$ has a dominated splitting.

Finally, the set
$$\mathcal{U}_1 = \bigcup_{f \in \mathcal{R}_0} \mathcal{U}(f)$$
is an open and dense subset of $\mathcal{U}$ satisfying the thesis of Lemma 2.0.2.

## 3. Proof of Theorem B

Recall from the introduction that $\Lambda$ has a dominated splitting $T_{/\Lambda}M = E \oplus F$ if there exist $C > 0$ and $0 < \lambda < 1$ such that
$$\|Df^n_{/E(x)}\|\|Df^{-n}_{/F(f^n(x))}\| < C\lambda^n$$
for all $x \in \Lambda$ and $n \geq 0$. Notice that the splitting must be continuous. We shall assume in this section that the constant $C$ above is equal to 1. If this is not the case, we can replace $f$ by a power of $f$: if Theorem B is true for a power of $f$ it is also true for $f$.

We shall reduce the proof of Theorem B to proving the following theorem.

THEOREM 3.1. *Let $f \in \mathrm{Diff}^2(M^2)$ and let $\Lambda$ be a compact invariant set having a dominated splitting $T_{/\Lambda}M = E \oplus F$. Also the periodic points of $f$ in $\Lambda$ are hyperbolic of saddle type. Then one of the following statements holds*:

1. $\Lambda$ *is a hyperbolic set*;

2. *There exists a simple closed curve* $\mathcal{C} \subset \Lambda$ *which is invariant under* $f^m$ *for some $m$ and is normally hyperbolic. Moreover $f^m : \mathcal{C} \to \mathcal{C}$ is conjugated to an irrational rotation.*

Before proving this last theorem, we show how it implies Theorem B. Let us assume that this last theorem is true. Assume that $\Lambda$ is a compact invariant set having a dominated splitting $T_{/\Lambda}M = E \oplus F$ and such that the periodic points of $f$ in $\Lambda$ are hyperbolic of saddle type. We will show that in this case, the number of periodic simple closed curves which are normally hyperbolic and conjugated to an irrational rotation contained in $\Lambda$ is finite. This will imply Theorem B.

Arguing by contradiction, assume that there exists an infinite number of such curves; that is, there exists a sequence $\mathcal{C}_n$ of simple closed curves invariant under $f^{m_n}$ for some $m_n$, normally hyperbolic, and such that $f^{m_n} : \mathcal{C}_n \to \mathcal{C}_n$ is conjugated to an irrational rotation. We can assume, without loss of generality, that these curves are normally attractive; i.e., for $x \in \mathcal{C}_n$, $T_x \mathcal{C}_n = F(x)$. The following lemma says that the diameter of these curves is bounded away from zero.

LEMMA 3.0.1. *There exist $\eta > 0$ such that $\mathrm{diam}(\mathcal{C}_n) > \eta$ for all $n$.*

*Proof.* Assume that the lemma is false. Then (taking a subsequence if necessary) we have

$$\mathrm{diam}(\mathcal{C}_n) \to_n 0.$$

Take a point $x_n \in \mathcal{C}_n$. Then there exist a subsequence $n_k$ such that $x_{n_k} \to x$ for some point $x \in \Lambda$. As $\mathrm{diam}(\mathcal{C}_n) \to_n 0$ we conclude that

$$\mathcal{C}_{n_k} \to x,$$

meaning that every sequence $y_{n_k} \in \mathcal{C}_{n_k}$ converges to the same point $x$. Since the dominated splitting is continuous, this implies that

$$F(y_{n_k}) \to F(x)$$

for every sequence $y_{n_k} \in \mathcal{C}_{n_k}$. But, for every $n_k$, there are points $y_{n_k}, z_{n_k} \in \mathcal{C}_{n_k}$ such that (trivializing the tangent bundle on a neighborhood of $x$) $F(y_{n_k})$ and $F(z_{n_k})$ are orthogonal. This is a contradiction because $F(y_{n_k}) \to F(x)$ and $F(z_{n_k}) \to F(x)$. □

We will get a contradiction of the existence of infinitely many such curves having diameters bounded away from zero. The idea is the following: since $f^{m_n} : \mathcal{C}_n \to \mathcal{C}_n$ is conjugated to an irrational rotation, they support only one invariant measure; thus the fiber $E$ is contractive. Then the basin of attraction of $\mathcal{C}_n$ has uniform size, and so, as there are infinitely many such curves, we will get an intersection of two different basins, which is a contradiction.

For this purpose we need the following lemma, which also will be used several times.

LEMMA 3.0.2. *Given a diffeomorphism $f$ and $0 < \gamma_1 < \gamma_2 < 1$, there exist a positive integer $N = N(\gamma_1, \gamma_2, f)$ and $c = c(\gamma_1, \gamma_2, f) > 0$ with the following property: For $x \in M$, a subspace $S \subset T_x M$ such that for some $n \geq N$ (with $S_i = Df^i(S)$),*

$$\prod_{i=0}^{n} \|Df_{/S_i}\| \leq \gamma_1^n,$$

*there exist* $0 \leq n_1 < n_2 < \cdots < n_l \leq n$ *such that*

$$\prod_{i=n_r}^{j} \|Df_{/S_i}\| \leq \gamma_2^{j-n_r}; \ r = 1, \ldots, l; \ n_r \leq j \leq n.$$

*Moreover,* $l \geq cn$.

We shall not prove this lemma because it is an immediate reformulation of a result by Pliss [Pl] (see also [M4, p. 276]).

Remember that in the case $\dim S = 1$, we have

$$\|Df_{/S}^n\| = \prod_{i=0}^{n-1} \|Df_{/S_i}\|.$$

We remark also that the dependence of $N$ and $c$ on $f$ is only in $\sup\{\|Df\|\}$. Modifying a little the constants we shall apply the lemma not only to $f$ but also to $f^{-1}$.

COROLLARY 3.1. *Given* $0 < \gamma_1 < \gamma_2 < 1$ *and* $x \in M$ *and a subspace* $S \subset T_xM$ *such that, for some* $m$,

$$\prod_{i=0}^{n} \|Df_{/S_i}\| \leq \gamma_1^n, \ \text{for all } n \geq m \ ,$$

*there exists a sequence* $0 \leq n_1 < n_2 < \cdots$ *such that*

$$\prod_{i=n_r}^{j} \|Df_{/S_i}\| \leq \gamma_2^{j-n_r}; \ \text{for all } j \geq n_r, r = 1, 2 \ldots \ .$$

We shall state the existence of locally invariant manifolds tangent to the $E$-direction and to the $F$-direction. For our immediate purpose we do not need these manifolds to be of class $C^2$, but in the followings sections this will play a fundamental role. To prove that they are of class $C^2$ we need the following lemma:

LEMMA 3.0.3. *Let* $\Lambda$ *be as in the statement of Theorem* B. *Then there exist a constant* $C > 0$ *and* $0 < \sigma < 1$ *such that for every* $x \in \Lambda$ *and for all positive integers* $n$ *the following hold*:

1. $\|Df^n_{/E(x)}\| \|Df^{-n}_{/F(f^n(x))}\|^2 < C\sigma^n$,

2. $\|Df^n_{/E(x)}\|^2 \|Df^{-n}_{/F(f^n(x))}\| < C\sigma^n$.

*Proof.* We shall prove only the first property of the thesis, the second one being analogous. First, it is enough to show the existence of a positive integer $m_1$ such that for every $x \in \Lambda$, there exists $1 \leq m \leq m_1$ such that

$$\|Df^m_{/E(x)}\| \|Df^{-m}_{/F(f^m(x))}\|^2 < \frac{1}{2}.$$

Now, arguing by contradiction, assume this is not true. Then, for each positive integer $n$ there exists $x_n \in \Lambda$ such that

$$\|Df^j_{/E(x_n)}\|\|Df^{-j}_{/F(f^j(x_n))}\|^2 \geq \frac{1}{2}$$

for all $0 \leq j \leq n$. We can assume, without loss of generality, that $x_n \to x$ for some $x \in \Lambda$. For this $x$ we get

$$\|Df^j_{/E(x)}\|\|Df^{-j}_{/F(f^j(x))}\|^2 \geq \frac{1}{2}$$

for all positive $j$. By the dominated splitting,

$$\|Df^j_{/E(x)}\|\|Df^{-j}_{/F(f^j(x))}\| < \lambda^j.$$

Then

$$\frac{1}{2} \leq \|Df^j_{/E(x)}\|\|Df^{-j}_{/F(f^j(x))}\|^2 < \|Df^{-j}_{/F(f^j(x))}\|\lambda^j.$$

Since $F$ is one dimensional, $\|Df^{-j}_{/F(f^j(x))}\| = \|Df^j_{/F(x)}\|^{-1}$ and so $\|Df^j_{/F(x)}\| < 2\lambda^j$. Again, by the dominated splitting,

$$\|Df^j_{/E(x)}\| < \lambda^j \|Df^j_{/F(x)}\| < 2\lambda^{2j}.$$

On the other hand, by the dominated splitting, angle$(E, F) > \gamma > 0$ for every point in $\Lambda$, and so there exists a positive constant $K$ such that

$$\|Df^n(z)\| \leq K \sup\{\|Df^n_{/E(z)}\|, \|Df^n_{/F(z)}\|\} \leq K\|Df^n_{/F(z)}\|$$

for all $z \in \Lambda$ and for all positive integer $n$, where the last inequality follows by the dominated splitting again. In particular, for every $q$,

$$\prod_{j=0}^{n} \|Df^q(f^{qj}(z))\| \leq \prod_{j=0}^{n} K\|Df^q_{/F(f^{qj}(z))}\|.$$

Take $\sigma_0$, $\lambda < \sigma_0 < 1$ and $q$ such that $2K\lambda^q < \sigma_0$. Then, for the point $x$,

$$\prod_{j=0}^{n} \|Df^q(f^{qj}(x))\| \leq K^n 2\lambda^{qn} = 2(K\lambda^q)^n \leq \sigma_0^n$$

for all $n \geq 0$. Let $g = f^q$. Thus

$$\prod_{j=0}^{n} \|Dg(g^j(x))\| \leq \sigma_0^n \text{ for all } n \geq 0.$$

Consider $0 < \lambda < \sigma_0 < \sigma_1 < \sigma_2 < 1$. Now, by Corollary 3.1, there exists a sequence of integers $n_k \to \infty$ such that, for any $k$ and for every positive $n$,

$$\|Dg^n(g^{n_k}(x))\| \leq \prod_{j=0}^{n-1} \|Dg(g^j(g^{n_k}(x)))\| < \sigma_1^n.$$

Thus, it can be proved that there exists $\eta > 0$, independent of $k$, such that for every $y, z \in B_\eta(g^{n_k}(x))$ we have $\mathrm{dist}(g^j(y), g^j(z)) \leq \sigma_2{}^j \mathrm{dist}(y, z)$ for every $j$. Let $j_0$ be such that for every $j \geq j_0$ we get $\sigma_2{}^j < \frac{\eta}{4}$. Now, take $n_i$ and $n_l$ such that $n_l - n_i > j_0$ and $\mathrm{dist}(g^{n_l}(x), g^{n_i}(x)) < \frac{\eta}{4}$. When we set $r = n_l - n_i$, it follows that $g^r(B_\eta(g^{n_i}(x))) \subset B_\eta(g^{n_i}(x))$ and also $g_{/B_\eta(g^{n_i}(x))}$ is a contraction. Then there is a point $p \in B_\eta(g^{n_i}(x))$ which is fixed under $g^r$, such that for every $z \in B_\eta(g^{n_i}(x))$ we have $g^{rn}(z) \to p$. Since $g = f^q$ we conclude that $p$ is an attracting fixed point under $f^{qr}$. Therefore $p$ is a sink, attracting the point $z = g^{n_i}(x)$ which is in $\Lambda$. Hence, $p \in \Lambda$, which is a contradiction (remember that we do not have sinks or repellors in $\Lambda$). $\square$

Let $I_1 = (-1, 1)$ and $I_\varepsilon = (-\varepsilon, \varepsilon)$, and denote by $\mathrm{Emb}^2(I_1, M)$ the set of $C^2$-embeddings of $I_1$ on $M$.

Recall by [HPS] that the dominated splitting and the previous lemma imply the next:

LEMMA 3.0.4. *There exist two continuous functions* $\phi^{cs} : \Lambda \to \mathrm{Emb}^2(I_1, M)$ *and* $\phi^{cu} : \Lambda \to \mathrm{Emb}^2(I_1, M)$ *such that, with* $W^{cs}_\varepsilon(x) = \phi^{cs}(x)I_\varepsilon$ *and* $W^{cu}_\varepsilon(x) = \phi^{cu}(x)I_\varepsilon$, *the following properties hold*:

a) $T_x W^{cs}_\varepsilon(x) = E(x)$ *and* $T_x W^{cu}_\varepsilon(x) = F(x)$,

b) *for all* $0 < \varepsilon_1 < 1$ *there exists* $\varepsilon_2$ *such that*

$$f(W^{cs}_{\varepsilon_2}(x)) \subset W^{cs}_{\varepsilon_1}(f(x))$$

*and*

$$f^{-1}(W^{cu}_{\varepsilon_2}(x)) \subset W^{cu}_{\varepsilon_1}(f^{-1}(x)).$$

We shall call the manifold $W^{cs}$ the (local) center stable manifold and $W^{cu}$ the (local) center unstable manifold. Observe that property b) means that $f(W^{cs}_\varepsilon(x))$ contains a neighborhood of $f(x)$ in $W^{cs}_\varepsilon(f(x))$ and $f^{-1}(W^{cu}_\varepsilon(x))$ contains a neighborhood of $f^{-1}(x)$ in $W^{cu}_\varepsilon(f^{-1}(x))$. In particular we have:

COROLLARY 3.2. *Given $\varepsilon$, there exists a number $\delta$ with the following properties*:

1. *If* $y \in W^{cs}_\varepsilon(x)$ *and* $\mathrm{dist}(f^j(x), f^j(y)) \leq \delta$ *for* $0 \leq j \leq n$ *then* $f^j(y) \in W^{cs}_\varepsilon(f^j(x))$ *for* $0 \leq j \leq n$.

2. *If* $y \in W^{cu}_\varepsilon(x)$ *and* $\mathrm{dist}(f^{-j}(x), f^{-j}(y)) \leq \delta$ *for* $0 \leq j \leq n$ *then* $f^{-j}(y) \in W^{cu}_\varepsilon(f^{-j}(x))$ *for* $0 \leq j \leq n$.

An open interval $I \subset M$ will mean for us an embedded one of the real line (or the open unit interval) in $M$. We denote by $\ell(I)$ its *length*.

COROLLARY 3.3. *Let $x \in \Lambda$ be such that, for some $0 < \gamma < 1$,*

$$\|Df^n_{/E(x)}\| \leq \gamma^n, \text{ for all } n \geq 0.$$

*Then, there exists $\varepsilon > 0$ such that*

$$\ell(f^n(W^{cs}_\varepsilon(x))) \to_n 0;$$

*i.e., the central stable manifold of size $\varepsilon$ is in fact a stable manifold.*

*Proof.* Take $\varepsilon_1 < 1$, and take $\delta = \delta(\varepsilon_1)$ as in the preceding corollary and small enough such that, if we define $\tilde{E}(z) = T_z W^{cs}_{\varepsilon_1}(y)$ for all $z \in W^{cs}_{\varepsilon_1}(y)$, $y \in \Lambda$, then,

$$\frac{\|Df_{/\tilde{E}(z_1)}\|}{\|Df_{/\tilde{E}(z_2)}\|} \leq 1 + c$$

for all $z_1, z_2$ such that $\text{dist}(z_1, z_2) < \delta$ where $(1+c)\gamma = \gamma_1 < 1$.

Let $\varepsilon$ be small enough such that $\ell(W^{cs}_\varepsilon(y)) \leq \delta$, for all $y \in \Lambda$. Now, it is not difficult to see that, for all $n \geq 0$,

$$\ell(f^n(W^{cs}_\varepsilon(x))) \leq \gamma_1^n \ell(W^{cs}_\varepsilon(x))$$

and

$$f^n(W^{cs}_\varepsilon(x)) \subset W^{cs}_{\varepsilon_1}(f^n(x)).$$

This implies the thesis.  □

In conclusion, we show how Theorem B follows from Theorem 3.1. Recall the assumption that there exists a sequence $\mathcal{C}_n$ of simple closed curves invariant under $f^{m_n}$ for some normally contractive $m_n$, such that $f^{m_n} : \mathcal{C}_n \to \mathcal{C}_n$ is conjugated to an irrational rotation. Moreover, we just proved that they have a diameter bounded away from zero.

Let $\lambda < \gamma_1 < \gamma_2 < 1$ and $c > 0$ such that $(1+c)\lambda < \gamma_1$. Since these curves are conjugated to an irrational rotation it is not difficult to see that for every $\mathcal{C}_n$ there exist points $x_n$ and $k_n$ such that

$$\|Df^j_{/F(x_n)}\| \leq (1+c)^j, \text{ for all } j \geq k_n.$$

Thus, by the dominated splitting,

$$\|Df^j_{/E(x_n)}\| \leq (1+c)^j \lambda^j < \gamma_1^j, \text{ for all } j \geq k_n$$

and so, by Corollary 3.1, there exists $j_n$ (we only need one, not a sequence) such that, setting $y_n = f^{j_n}(x_n)$, we have

$$\|Df^j_{/E(y_n)}\| \leq \gamma_2^j, \text{ for all } j \geq 0.$$

Now, by Corollary 3.3 there exists $\varepsilon$ such that, for every $n$,

$$\ell(f^j(W^{cs}_\varepsilon(y_n))) \to_j 0.$$

Then, since the diameters of the curves are bounded away from zero, we can find $n_1 \neq n_2$ big enough, such that $W^{cs}_\varepsilon(y_{n_1})$ intersects the orbit of $\mathcal{C}_{n_2}$ (remember that $F$ is the tangent direction of $\mathcal{C}_n$, for all $n$). This is a contradiction, because the point of intersection remains in the orbit of $\mathcal{C}_{n_2}$ and also is asymptotic to the orbit of $\mathcal{C}_{n_1}$.

Thus, there exist only finitely many periodic simple closed curves normally hyperbolic and conjugated to an irrational rotation contained in $\Lambda$. Let $\Lambda_2$ be the union of such curves. Since these curves are isolated, we conclude that $\Lambda_1 = \Lambda - \Lambda_2$ is a compact invariant set. Now applying Theorem 3.1 we conclude that $\Lambda_1$ is a hyperbolic set. Hence Theorem 3.1 implies Theorem B.

The proof of Theorem 3.1 will be given in the following section.

3.1. *Proof of Theorem* 3.1. The first step in the proof is the following elementary lemma. The proof is left to the reader.

LEMMA 3.1.1. *Let $\Lambda_0$ be a compact invariant set having a dominated splitting $T_{/\Lambda_0}M = E \oplus F$. If $\|Df^n_{/E(x)}\| \to 0$ and $\|Df^{-n}_{/F(x)}\| \to 0$ as $n \to \infty$ for every $x \in \Lambda_0$ then $\Lambda_0$ is a hyperbolic set.*

Now, we will prove Theorem 3.1 based on the next lemma.

MAIN LEMMA. *Let $\Lambda_0$ be a nontrivial, transitive and compact, invariant set having a dominated splitting $T_{/\Lambda_0}M = E \oplus F$ and such that it is not a periodic simple closed curve normally hyperbolic $\mathcal{C}$ conjugated to an irrational rotation. Assume that every properly compact invariant subset of $\Lambda_0$ is hyperbolic. Then, $\Lambda_0$ is a hyperbolic set.*

The proof of this lemma will be given in subsection 3.7.

Let us prove how the Main Lemma implies Theorem 3.1. For that, it is necessary to show that either 1) or 2) of the statement holds. Assume that 2) does not hold; that is, there is not a periodic simple closed curve $\mathcal{C} \subset \Lambda$ normally hyperbolic conjugated to an irrational rotation.

Then we have to prove that $\Lambda$ is a hyperbolic set. Arguing by contradiction, we assume this is not true.

Let
$$H = \{\tilde{\Lambda} \subset \Lambda : \tilde{\Lambda} \text{ is a nonhyperbolic compact invariant set}\},$$
and we order $H$ by inclusion.

Let $H_\Gamma = \{\Lambda_\gamma : \gamma \in \Gamma\}$ be a totally ordered chain. Then $\Lambda_\infty = \cap_{\gamma \in \Gamma}\Lambda_\gamma$ is a compact invariant set. We claim that $\Lambda_\infty$ is also a nonhyperbolic set. If it were hyperbolic, then we could take a small neighborhood $U$ of $\Lambda_\infty$ in such a way that every compact invariant set in $U$ would be hyperbolic. By the definition of $\Lambda_\infty$, we would have, for some $\gamma \in \Gamma$, that $\Lambda_\gamma$ is contained in $U$, implying that it is a hyperbolic set, which is a contradiction.

Now, by Zorn's lemma there exists a nonhyperbolic compact invariant set $\Lambda_0 \subset \Lambda$ such that every properly compact invariant subset is hyperbolic. We claim that $\Lambda_0$ is a nontrivial transitive set, that is, a transitive set which is not a periodic orbit.

First of all, since the periodic points in $\Lambda$ are hyperbolic, $\Lambda_0$ cannot be a periodic orbit. Secondly, if for every point $x \in \Lambda_0$ we have that the $\alpha$-limit set $\alpha(x)$ is properly contained in $\Lambda_0$, then it will be hyperbolic. This implies that $\|Df^{-n}_{/F(x)}\| \to 0$ as $n \to \infty$ for every point $x$. In the same way, if the $\omega$-limit set of any point $x \in \Lambda_0$ is properly contained in $\Lambda_0$, we conclude that $\|Df^n_{/E(x)}\| \to 0$ as $n \to \infty$ for every point $x$. But this implies that $\Lambda_0$ (see Lemma 3.1.1) is hyperbolic which is a contradiction, proving that $\Lambda_0$ is a transitive set.

In conclusion, the set $\Lambda_0$ satisfies the hypothesis of the Main Lemma but as it was defined is not hyperbolic which yields a contradiction. Thus, Theorem 3.1 follows from the Main Lemma.

As we just said, the proof of the Main Lemma will be given in Section 3.7. Nevertheless, we give here the basic steps of it:

1. The central stable manifolds (which are of class $C^2$) have dynamics properties. In fact for every $x \in \Lambda_0$ there exists $\varepsilon(x)$ such that $W^{cs}_{\varepsilon(x)}(x)$ is a stable manifold of $x$, and $W^{cu}_{\varepsilon(x)}(x)$ an unstable manifold of $x$, meaning that $\ell(f^n(W^{cs}_{\varepsilon(x)}(x))) \to 0$ and $\ell(f^{-n}(W^{cu}_{\varepsilon(x)}(x))) \to 0$ as $n \to \infty$.

2. There exists an open set $B$ in $\Lambda_0$ such that, for every $x \in B$,
$$\sum_{n \geq 0} \ell(f^n(W^{cs}_{\varepsilon(x)}(x))) < \infty$$
and
$$\sum_{n \geq 0} \ell(f^{-n}(W^{cu}_{\varepsilon(x)}(x))) < \infty.$$

3. For every point $x \in \Lambda_0$,
$$\|Df^n_{/E(x)}\| \to 0$$
and
$$\|Df^{-n}_{/F(x)}\| \to 0$$
when $n \to \infty$.

3.2. *A Denjoy property.* The aim of this section is to prove that in a set with a dominated splitting, there is no wandering (in the future) interval transversal to the $E$-direction. We must state some basic properties on dominated splitting. First we introduce the notion of cone fields.

Let $\Lambda$ be a compact invariant set having a dominated splitting $T_{/\Lambda}M = E \oplus F$. For $0 < a < 1$ we define the cones:

$$C^{cu}_a(x) = \{w \in T_xM : w = v_E + v_F, v_E \in E(x), v_F \in F(x), \|v_E\| \leq a\|v_F\|\}$$

and

$$C^{cs}_a(x) = \{w \in T_xM : w = v_E + v_F, v_E \in E(x), v_F \in F(x), \|v_F\| \leq a\|v_E\|\}.$$

LEMMA 3.2.1. *Assume that $\Lambda$ is a compact invariant set with a dominated splitting $T_{/\Lambda}M = E \oplus F$ such that*

$$\|Df_{/E(x)}\|\|Df^{-1}_{/F(f(x))}\| < \lambda.$$

*Then, for any $0 < a \leq 1$ and $x \in \Lambda$, the following hold:*

$$Df_x.C_a^{cu}(x) \subset C_{\lambda a}^{cu}(f(x)),$$
$$Df_x^{-1}.C_a^{cs}(x) \subset C_{\lambda a}^{cs}(f^{-1}(x)).$$

*Conversely, if a compact invariant set $\Lambda$ has a continuous decomposition $T_\Lambda M = E_1 \oplus F_1$ (not necessarily invariant) and for some $0 < a \leq 1$ it is true that*

$$Df_x.C_a^{cu}(x) \subset C_{\lambda a}^{cu}(f(x))$$

*and*

$$Df_x^{-1}.C_a^{cs}(x) \subset C_{\lambda a}^{cs}(f^{-1}(x)),$$

*then $\Lambda$ has a dominated splitting (here the cones are defined relative to the decomposition $T_\Lambda M = E_1 \oplus F_1$).*

We shall not prove this lemma because it is similar to the hyperbolic case.

Now, let $\Lambda$ be a compact invariant set having a dominated splitting as in Lemma 3.2.1, and take cone fields $C_a^{cu}$, $C_a^{cs}$ (we shall refer to them as the central unstable and stable $a$-cone respectively). Then, there exists an admissible neighborhood $V_a(\Lambda)$; i.e., we can extend these cones to $V_a(\Lambda)$ having the same properties as for $\Lambda$ where it make sense.

*Definition* 1. An interval $I \subset V_a(\Lambda)$ is $a$-transversal to the $E$-direction if for every $x \in I$, we have that $T_xI$ is contained in the $C_a^{cu}(x)$. We say that it is transversal if it is $a$-transversal for some $0 < a \leq 1$.

*Remark* 3.1. Notice that there exists $\eta$ such that if $I$ is an interval transversal to the $E$-direction and $\ell(I) < \eta$ then there exists another interval $I_0$ transversal to the $E$-direction, containing $I$ and $\ell(I_0) \geq \eta$. In the sequel, the number $\delta > 0$ will always be less than $\eta$. Moreover, notice that central unstable manifolds are always transversal to the $E$-direction.

From now on, we fix an admissible neighborhood $V = V_1$. Take $U$ another neighborhood of $\Lambda$ such that $U \subset \text{Cl}(U) \subset V$ (where Cl denotes the closure). Denote by $\Lambda_1 = \cap_{n \in \mathbb{Z}} f^n(\text{Cl}(U))$ the maximal invariant set in $\text{Cl}(U)$ and by $\Lambda_1^+ = \cap_{n \geq 0} f^{-n}(\text{Cl}(U))$ the set of points which remains in $\text{Cl}(U)$ in the future. We remark that $\Lambda_1$ has a dominated splitting $T_{\Lambda_1} = E \oplus F$ since $V$ is admissible. Moreover, for every point $x \in \Lambda_1^+$ we have a uniquely determined $E$-direction. Let us make the following:

*Definition* 2. We say that an open $C^2$ interval $I$ in $M$ is a $\delta$-$E$-interval if the next two conditions hold:

1. $I \subset \Lambda_1^+$ and $\ell(f^n(I)) \leq \delta$ for all $n \geq 0$.

2. $f^n(I)$, $n \geq 0$, is always transversal to the $E$-direction.

Furthermore, if the supplementary condition:

3. There exists $\gamma > 0$ such that for every $x \in I$ and $n \geq 0$ $\|Df^n_{/E(x)}\| < \gamma^n$ is satisfied, we say that $I$ is a $(\delta, \gamma)$-$E$-interval.

In an analogous way we define the $\delta$-$F$-interval and the $(\delta, \gamma)$-$F$-interval.

LEMMA 3.2.2. *Let $0 < \lambda < \gamma$. Then, there exists $\delta_1 = \delta_1(\gamma) > 0$ such that if $I$ is a $\delta$-$E$-interval with $\delta \leq \delta_1$ then, for some $m \geq 0$, $f^m(I)$ is a $(\delta, \gamma)$-$E$-interval.*

*Proof.* Take $0 < \lambda < \lambda_0 < \lambda_1 < \gamma < \lambda_2 < \lambda_3 < 1$ such that $\lambda \lambda_2^{-1} < \lambda_0$ and take $c > 0$ such that $(1+c)\lambda_1 < \gamma$ and $(1-c)\lambda_2^{-1} > \lambda_3^{-1}$. Take now $\delta_1 > 0$ and $V_a(\Lambda_1)$, with $a$ sufficiently small, such that

$$\frac{\|Df_{/\tilde{F}(x)}\|}{\|Df_{/\tilde{F}(y)}\|} > 1 - c$$

and

$$\frac{\|Df_{/E(x)}\|}{\|Df_{/E(y)}\|} < 1 + c$$

for every $x, y$ in $V_a(\Lambda_1)$ with $\text{dist}(x, y) < \delta_1$, where $\tilde{F}$ is any direction in the central unstable $a$-cone. Let $I$ be a $\delta$-$E$-interval with $\delta \leq \delta_1$. Observe that to prove the lemma it is enough to find some $m$ such that $f^m(I)$ satisfies the third property, because the first two are already satisfied for every $m \geq 0$.

Since $I \subset \Lambda_1^+$ and the $w$-limit set of every point of $\Lambda_1^+$ is contained in $\Lambda_1$, we conclude that $f^n(I) \subset V_a(\Lambda_1)$ for all suficiently large $n$. Moreover, $T_x f^n(I)$ is contained in the central unstable $a$-cone. Thus, for our purpose, we can assume that these facts hold for any $n \geq 0$.

Let $\tilde{F}(x) = T_x f^n(I)$ for $x \in f^n(I)$. Since $\tilde{F}(x) \subset C_a^{cu}(x)$ and $a$ is small, it follows that we have the domination property

$$\|Df_{/E(x)}\| \|Df^{-1}_{/\tilde{F}(f(x))}\| < \lambda.$$

Take $x \in I$. We claim that for all arbitrary large $n$ we have $\|Df^n_{/\tilde{F}(x)}\| \leq \lambda_2^{-n}$. Otherwise, there exists a sequence $n_k \to \infty$ such that $\|Df^{n_k}_{/\tilde{F}(x)}\| > \lambda_2^{-n_k}$. Then, for every $y \in I$,

$$\|Df^{n_k}_{/\tilde{F}(y)}\| > ((1-c)\lambda_2^{-1})^{n_k} > \lambda_3^{-n_k}.$$

Taking $n_k$ such that $\lambda_3^{-n_k}\ell(I) > \delta$, we conclude

$$\ell(f^{n_k}(I)) > \lambda_3^{-n_k}\ell(I) > \delta,$$

which is a contradiction because $I$ is a $\delta$-$E$-interval. Thus, $\|Df^n_{/\tilde{F}(x)}\| \leq \lambda_2^{-n}$ for every large $n$. By the domination property we conclude

$$\|Df^n_{/E(x)}\| \leq (\lambda\lambda_2^{-1})^n \leq \lambda_0^n$$

for every large $n$. By Corollary 3.1 there exists a positive integer $m$ such that $\|Df^n_{/E(f^m(x))}\| \leq \lambda_1^n$ for every positive $n$. Finally, this implies that for every $y \in f^m(I)$ it is true that $\|Df^n_{/E(y)}\| \leq \gamma^n$ for every positive $n$. Hence $f^m(I)$ is a $(\delta, \gamma)$-$E$-interval. □

We say that an interval $J$ has a stable property if for every $x \in J$ there exists a local stable manifold $W^s_\varepsilon(x)$ of uniform size. We shall call an interval $I$ pre-periodic if there exists an interval $J$ with a stable property such that for some $m > 0$ $f^m(J) \subset J$ and for some $n \geq 0$ we have $f^n(I) \subset \cup_{x \in J} W^s_\varepsilon(x)$. It is not difficult to see that if $I$ is pre-periodic, then $\omega(x)$ is a periodic orbit for every $x \in I$. For an interval $I$ denote by $\omega(I) = \cup_{x \in I} \omega(x)$ the union of the $\omega$-limit set of the points in $I$.

The next result will play a central role. It establish a kind of Denjoy theorem.

PROPOSITION 3.1. *There exists $\delta_0$ such that if $I$ is a $\delta$-$E$-interval with $\delta \leq \delta_0$, then one of the following properties holds*:

1. $\omega(I)$ is a periodic simple closed curve $\mathcal{C}$ normally hyperbolic and $f^m_{/\mathcal{C}} : \mathcal{C} \to \mathcal{C}$ (where $m$ is the period of $\mathcal{C}$) is conjugated to an irrational rotation,

2. $\omega(I) \subset \operatorname{Per}(f_{/V})$ where $\operatorname{Per}(f_{/V})$ is the set of the periodic points of $f$ in $V$.

*Proof.* Take $\gamma = \lambda^{\frac{1}{2}}$, and $\lambda_2, \lambda_3$; $\lambda < \gamma < \lambda_2 < \lambda_3 < 1$. Also, let $c > 0$ be such that $(1+c)\lambda_2 < \lambda_3$.

Pick $\delta_2 > 0$ and $a$ small such that if $\operatorname{dist}(x,y) < \delta_2$, $x, y \in V_a(\Lambda_1)$ then

$$(1+c)^{-1} < \frac{\|Df_{/E(x)}\|}{\|Df_{/E(y)}\|} < 1 + c$$

and

$$(1+c)^{-1} < \frac{\|Df^{-1}_{/\tilde{F}(x)}\|}{\|Df^{-1}_{/\tilde{F}(y)}\|} < 1 + c$$

where $\tilde{F}$ is any direction in the central unstable $a$-cone. Take $\delta_1 = \delta_1(\gamma)$ from the previous lemma and let $\delta_0 \leq \min\{\delta_1, \delta_2\}$.

Let $I$ be a $\delta$-$E$-interval with $\delta \leq \delta_0$. By the previous lemma, for some $m \geq 0$, $f^m(I)$ is a $(\delta, \gamma)$-$E$-interval. It suffices to prove the proposition for $f^m(I)$. Thus, we may assume that $m = 0$; that is, $I$ is a $(\delta, \gamma)$-$E$-interval. Moreover, as in the previous lemma, if we set $\tilde{F}(x) = T_x f^n(I)$, $x \in f^n(I)$, the domination property can be assumed:

$$\|Df_{/E(x)}\|\|Df^{-1}_{/\tilde{F}(f(x))}\| < \lambda.$$

Take some $x_0 \in I$. Then, for every positive $n$ we have $\|Df^n_{/E(x_0)}\| < \gamma^n$. Take the sequence of all integers $n_i$, $n_i \to \infty$, such that

$$\|Df^j_{/E(f^{n_i}(x_0))}\| < \lambda_2^j \text{ for all } j \geq 0$$

(such a sequence exists by Corollary 3.1). This implies, by the way we choose $\delta_0$, that $f^{n_i}(I)$ is a $(\delta, \lambda_3)$-$E$-interval. Take now $I_0$, $I \subset I_0$, a maximal $(\delta, \lambda_3)$-$E$-interval. It follows that $f^{n_i}(I_0)$ is a $(\delta, \lambda_3)$-$E$-interval for every $n_i$. Now, for $i \geq 1$, take $I_{n_i}$ such that $f^{n_i - n_{i-1}}(I_{n_{i-1}}) \subset I_{n_i}$ and $I_{n_i}$ is a maximal $(\delta, \lambda_3)$-$E$-interval.

By Corollary 3.3 we have, for every $i$, that the interval $I_{n_i}$ has the stable property. Moreover, $W^s_\varepsilon(x)$ has uniform length for every $x \in I_{n_i}$ and for every $i$. Furthermore, it can be proved that there exists a constant $K$ such that, if the box $W^s_\varepsilon(I_{n_i}) = \cup_{x \in I_{n_i}} W^s_\varepsilon(x)$ is defined, we have $K\text{vol}(W^s_\varepsilon(I_{n_i})) \geq \ell(I_{n_i})$, where $K$ is independent of $i$ (see also §§3.4 and the proof of Lemma 3.7.3).

Next, assume that there exist $n_i < n_j$ such that

$$W^s_\varepsilon(I_{n_i}) \cap W^s_\varepsilon(f^{n_j - n_i}(I_{n_i})) \neq \emptyset.$$

Let $m = n_j - n_i$. If $\ell(f^{km}(I_{n_i})) \to 0$ as $k \to \infty$, then $\omega(I_{n_i})$ consists of a periodic orbit. Indeed, if $\ell(f^{km}(I_{n_i})) \to 0$, then $\ell(f^k(I_{n_i})) \to 0$ as $k \to \infty$. Let $p$ be an accumulation point of $f^k(I_{n_i})$; that is, $f^{k_j}(I_{n_i}) \to p$ for some $k_j \to \infty$, and so, $f^{k_j + m}(I_{n_i}) \to f^m(p)$. But by the property we are assuming, i.e., $W^s_\varepsilon(I_{n_i}) \cap W^s_\varepsilon(f^{n_j - n_i}(I_{n_i})) \neq \emptyset$, we have $f^{k_j + m}(I_{n_i}) \to p$, implying that $p$ is a periodic point. Thus, for any $x \in I_{n_i}$ we have that $\omega(x)$ consists only of periodic orbits, and so $\omega(x)$ is a single periodic orbit $p$. Since $\ell(f^k(I_{n_i})) \to 0$ we conclude that $\omega(I_{n_i})$ is the orbit of the periodic point $p$. By the way we choose $I_{n_i}$, we have $f^{n_i}(I) \subset I_{n_i}$ and so $\omega(I)$ consists of a periodic orbit, as the thesis of the proposition requires.

On the other hand, if $\ell(f^{km}(I_{n_i}))$ does not go to zero, we take a sequence $k_j$ such that $f^{k_j m}(I_{n_i}) \to L$ for some interval $L \subset \Lambda_1$ (which is at least $C^1$, and has $F$ as its tangent direction). Now $f^{(k_j + 1)m}(I_{n_i}) \to L'$ and $f^m(L) = L'$. Moreover, by the property we are dealing with, $L \cup L'$ is an interval (with $F$ as its tangent direction). Let

$$J = \bigcup_{n \geq 0} f^{nm}(L).$$

We claim that there are only two possibilities: either $J$ is an interval or a simple closed curve. To prove this, notice that $f^{nm}(L)$ is a $\delta$-$E$-interval for any $n \geq 0$. In particular, for any $x \in J$ there exists $\varepsilon(x)$ such that $W^{cs}_{\varepsilon(x)}(x)$ is a stable manifold for $x$, and so

$$W(J) = \bigcup_{x \in J} W^{cs}_{\varepsilon(x)}(x)$$

is a neighborhood of $J$.

We only have to show that, given $x \in J$, there exists a neighborhood $U(x)$ such that $U(x) \cap J$ is an interval. This implies that $J$ is a simple closed curve or an interval. Thus, take $x \in J$, in particular $x \in f^{n_1 m}(L)$. Take $U$ an open interval, $x \in U \subset f^{n_1 m}(L)$ and let $U(x)$ be a neighborhood of $x$ such that $U(x) \subset W(J)$ and $U(x) \cap L_1 \subset U$ where $L_1$ is any interval containing $f^{n_1 m}(L)$, transversal to the $E$-direction and $\ell(L_1) \leq 2\delta_0$ (this is always possible if $\delta_0$ is small). Now let $y \in J \cap U(x)$. We have to prove that $y \in U$. There is $n_2$ such that $y \in f^{n_2 m}(L)$. Since

$$f^{n_1 m}(L) = \lim_j f^{k_j m + n_1 m}(I_{n_i}),$$
$$f^{n_2 m}(L) = \lim_j f^{k_j m + n_2 m}(I_{n_i})$$

and both have nonempty intersection with $U(x)$, we conclude that for some $j$, $f^{k_j m + n_1 m}(I_{n_i})$ and $f^{k_j m + n_2 m}(I_{n_i})$ are linked by a local stable manifold. Hence $f^{n_1 m}(L) \cup f^{n_2 m}(L)$ is an interval $L_1$ transversal to the $E$-direction with $\ell(L_1) \leq 2\delta_0$. Therefore $y \in U(x) \cap L_1 \subset U$ as we wish, completing the proof that $J$ is an interval or a simple closed curve.

In case $J$ is an interval, since $f^m(J) \subset J$, it follows that $I_{n_i}$ (and so $I$) is a pre-periodic interval and so, for any $x \in I$, $\omega(x)$ is a periodic orbit, which completes the proof in this case. On the other hand, if $J$ is a simple closed curve, which is of class $C^2$ because it is normally hyperbolic (attractive), then we have two possibilities. If $f^m_{/J} : J \to J$ has a rational rotation number, then we can see that the $\omega(I_{n_i})$ consist of a union of periodic points, and the same happens to $I$. If $f^m_{/J} : J \to J$ has an irrational rotation number, then it is conjugated to an irrational rotation, and denoting $\mathcal{C} = J$, we have that $\omega(I)$ is as in the first conclusion of the proposition.

Thus, we have proved the lemma in case there exist some $n_i < n_j$, such that $W^s_\varepsilon(I_{n_i}) \cap W^s_\varepsilon(f^{n_j - n_i}(I_{n_i})) \neq \emptyset$. Now, assume this is not true. We claim that

$$\sum_{k \geq 0} \ell(f^k(I_{n_i})) < \infty$$

for every $i$. Let us prove the claim.

Due to the fact that $f^{n_i-n_{i-1}}(I_{n_{i-1}}) \subset I_{n_i}$,

$$W^s_\varepsilon(I_{n_i}) \cap W^s_\varepsilon(f^{n_j-n_i}(I_{n_i})) = \emptyset$$

implies

$$W^s_\varepsilon(f^{n_k-n_i}(I_{n_i})) \cap W^s_\varepsilon(f^{n_j-n_i}(I_{n_i})) = \emptyset$$

for every $n_k > n_i, n_j > n_i$ and for every $i$.

Let us fix some $i$ and prove our claim. For simplicity we assume that $i = 0$. Let $N = N(\gamma, \lambda_2)$ be as in Lemma 3.0.2, and $m_j = n_{j+1} - n_j$. Then $m_j \leq N$ or $m_j > N$.

Assume that $m_j > N$. Then,

$$\|Df^n_{/E(f^{n_{j+1}-n}(x_0))}\| \geq \gamma^n$$

for $N \leq n \leq m_j$. Otherwise, if for some $N \leq n \leq m_j$,

$$\|Df^n_{/E(f^{n_{j+1}-n}(x_0))}\| < \gamma^n,$$

then there exists some $0 < \tilde{n} < n$ such that

$$\|Df^j_{/E(f^{n_{j+1}-\tilde{n}}(x_0))}\| < \lambda_2^j$$

for $0 \leq j \leq \tilde{n}$ by Lemma 3.0.2. But since $f^{\tilde{n}}(f^{n_{j+1}-\tilde{n}}(x_0)) = f^{n_{j+1}}(x_0)$ we conclude

$$\|Df^j_{/E(f^{n_{j+1}-\tilde{n}}(x_0))}\| < \lambda_2^j$$

for all positive $j$ which is a contradiction to the way we chose the sequence $n_i$, because $n_j < n_{j+1} - \tilde{n} < n_{j+1}$.

Now, if for $N \leq n \leq m_j$

$$\|Df^n_{/E(f^{n_{j+1}-n}(x_0))}\| \geq \gamma^n$$

then, by the domination property,

$$\|Df^{-n}_{/\tilde{F}(f^{n_{j+1}}(x_0))}\| \leq (\lambda\gamma^{-1})^n = \gamma^n$$

for all $N \leq n \leq m_j$. Using Lemma 3.0.2 again, we conclude that there exists some $\tilde{n}_j, n_{j+1} - N < \tilde{n}_j \leq n_{j+1}$ such that

$$\|Df^{-j}_{/\tilde{F}(f^{\tilde{n}_j}(x_0))}\| < \lambda_2^j$$

for $0 \leq j \leq \tilde{n}_j - n_j$. By the way we choose $\delta_0$ we get, for all $y \in f^{\tilde{n}_j}(I_0)$, $\|Df^{-j}_{/\tilde{F}(y)}\| < \lambda_3^j$ for $0 \leq j \leq \tilde{n}_j - n_j$. But then

$$\sum_{j=0}^{\tilde{n}_j-n_j} \ell(f^{-j}(f^{\tilde{n}_j}(I_0))) \leq \sum_{j=0}^{\tilde{n}_j-n_j} \lambda_3^j \ell(f^{\tilde{n}_j}(I_0)) \leq \frac{1}{1-\lambda_3}\ell(f^{\tilde{n}_j}(I_0)).$$

Let $K_1 = \sup\{\|Df^{-n}(x)\| : x \in M, 0 \le n \le N\}$ and $K_2 = \frac{1}{1-\lambda_3}K_1 + NK_1$. Now we get

$$\begin{aligned}
\sum_{k \ge 0} \ell(f^k(I_0)) &= \sum_j \sum_{n_j}^{n_{j+1}-1} \ell(f^k(I_0)) \\
&= \sum_{j, m_j > N} \sum_{n_j}^{n_{j+1}-1} \ell(f^k(I_0)) + \sum_{j, m_j \le N} \sum_{n_j}^{n_{j+1}-1} \ell(f^k(I_0)) \\
&\le \sum_{j, m_j > N} (\frac{1}{1-\lambda_3} K_1 \ell(f^{n_{j+1}}(I_0)) + NK_1 \ell(f^{n_{j+1}}(I_0))) \\
&\quad + \sum_{j, m_j \le N} NK_1 \ell(f^{n_{j+1}}(I_0)) \\
&\le K_2 \sum_j \ell(f^{n_j}(I_0)) \le KK_2 \sum_j \mathrm{vol}(W^s_\varepsilon(f^{n_j}(I_0))) \\
&\le KK_2 \mathrm{vol}(M) < \infty.
\end{aligned}$$

This proves our claim. In particular, $\ell(f^k(I_{n_i})) \to 0$ as $k \to \infty$ for all $i$.

Let us continue with the proof of the proposition. If

$$W^s_\varepsilon(I_{n_r}) \cap W^s_\varepsilon(I_{n_j}) \ne \emptyset$$

for some $n_r \ne n_j$ then $\omega(I)$ will be a periodic orbit. In fact, assuming $n_r < n_j$, using the fact that $f^{n_j - n_r}(I_{n_r}) \subset I_{n_j}$ and also that $\ell(f^k(I_{n_i})) \to 0$ for every $i$, we can see (as we did before) that if $p$ is an accumulation point of the sequence $\{f^k(I_{n_r})\}_k$ then it will be a periodic point with period $n_j - n_r$. Hence, $\omega(I)$ is a periodic orbit.

On the other hand, assume

$$W^s_\varepsilon(I_{n_r}) \cap W^s_\varepsilon(I_{n_j}) = \emptyset$$

for every $r, j$. This implies in particular that $\ell(I_{n_i}) \to 0$ as $i \to \infty$. Now, using Schwartz's method in the proof of the Denjoy theorem ([Sch]), we get a contradiction to the maximality of $I_{n_i}$ or for every $i$, there exists $k_i$ such that

$$\ell(f^{k_i}(I_{n_i})) = \delta.$$

We claim in the latter, that there exists $n_j > n_i$ such that $k_i = n_j - n_i$. This implies, in particular, that $\ell(I_{n_j}) = \delta$ which is a contradiction to the fact that $\ell(I_{n_i}) \to 0$ as $i \to \infty$.

So, to finish the proof of the proposition it is sufficient to prove our last claim. For this purpose, we only have to show that

$$\|Df^j_{/E(f^{n_i + k_i}(x_0))}\| \le \lambda_2^j$$

for every $j \ge 0$. Arguing by contradiction, we assume that for some $j_0 > 0$,

$$\|Df^{j_0}_{/E(f^{n_i + k_i}(x_0))}\| > \lambda_2^{j_0}.$$

This implies, by the domination property, that

$$\|Df^{-j_0}_{/\tilde{F}(f^{n_i+k_i+j_0}(x_0))}\| \leq (\lambda\lambda_2^{-1})^{j_0} \leq \gamma^{j_0} \leq \lambda_2^{j_0}.$$

Therefore, for every $y \in f^{j_0}(f^{k_i}(I_{n_i}))$,

$$\|Df^{-j_0}_{/\tilde{F}(y)}\| \leq (\lambda^2(1+c))^{j_0} \leq \lambda_3^{j_0}$$

and so

$$\|Df^{j_0}_{/F(y)}\| \geq \lambda_3^{-j_0}$$

for every $y \in f^{k_i}(I_{n_i})$. Then

$$\delta \geq \ell(f^{j_0}(f^{k_i}(I_{n_i}))) \geq \lambda_3^{-j_0}\ell(f^{k_i}(I_{n_i})) = \lambda_3^{-j_0}\delta > \delta,$$

which is a contradiction. This proves our last claim and completes the proof of Proposition 3.1. □

COROLLARY 3.4. *Assume that there is not a simple closed curve $\mathcal{C}$, which is $f^m$-invariant for some $m$, contained in $\Lambda$. Then there is an admissible neighborhood $V$ of $\Lambda$ such that any $\delta$-E-interval $I$, with $\delta$ sufficiently small, has $\omega(I) \subset \text{Per}(f_{/V})$.*

*Proof.* Let $V_0$ be any admissible neighborhood of $\Lambda$. Take another neighborhood $U$, such that $U \subset \text{Cl}(U) \subset V$ and let $M(U)$ be the maximal invariant set in U. If there is no curve in $U$ as in the statement, then by the previous proposition we conclude the proof. If there are such curves, since $M(U)$ has a dominated splitting, there can be only a finite number of them (see the beginning of the proof of Theorem B). Since no one is in $\Lambda$, another neighborhood $V$ can be taken in which there is no such curve. This $V$ satisfies the thesis. □

*Remark* 3.2. All the results proved in this section have similar versions for $\delta$-F-intervals when we replace $f$ by $f^{-1}$ and $E$ by $F$; the proofs are analogous.

3.3. *Dynamical properties of $W^{cu}$, $W^{cs}$.* From now on, we shall assume that $\Lambda$ is a compact invariant set with a dominated splitting such that there is no simple closed curve $\mathcal{C} \subset \Lambda$, $f^m$-invariant for some $m$, normally hyperbolic and conjugated to an irrational rotation.

Recall, by the previous section, that in this case there exists $\delta_0$ such that, for any $\delta$-E-interval $I$ with $\delta \leq \delta_0$, $\omega(I) \subset \text{Per}(f_{/V})$, where $V$ is an appropriate admissible neighborhood. We shall assume that $\delta_0 \leq \delta_3$ where $\delta_3$ is as in the following:

LEMMA 3.3.1. *There exists $\delta_3 > 0$ such that if $p \in \Lambda$ is a periodic point, and one component of $W^u(p) - \{p\}$ has length less than $\delta_3$, then the other endpoint of this component is a periodic point which is not hyperbolic of saddle type; more than that it is a sink or a nonhyperbolic periodic point.*

*Proof.* Let $V$ be an admissible neighborhood of $\Lambda$. Take $\delta_3 > 0$ such that

$$\{x : \text{dist}(x, \Lambda) \leq \delta_3\} \subset V.$$

Let $p \in \Lambda$ be a periodic point such that one component of $W^u(p) - \{p\}$ (say $W$) has length less than $\delta_3$. We may assume that $p$ is a fixed point. Then, $W$ is invariant and so $W \subset \Lambda_1$, the maximal invariant subset in $V$, which has a dominated splitting $T_{\Lambda_1} M = E \oplus F$. Since $T_p W^u(p) = F(p)$ and the splitting is continuous, we conclude that $T_x W = F(x)$ for any $x \in W$. Hence, $W$ is an interval transversal to the $E$-direction. Let $q$ be the other endpoint of $W$. Thus, $q$ is a fixed point attracting all the points in $W$. Therefore,

$$\|Df_{/F(q)}\| \leq 1.$$

By the domination,

$$\|Df_{/E(q)}\| \leq \lambda < 1,$$

which concludes the proof of the lemma. $\square$

Dynamical properties on the central unstable and stable manifolds are the main consequence of the preceding section:

LEMMA 3.3.2. *For all $\varepsilon < \delta_0$ there exists $\gamma = \gamma(\varepsilon)$ such that:*

1. *$f^{-n}(W_\gamma^{cu}(x)) \subset W_\varepsilon^{cu}(f^{-n}(x))$ and $f^n(W_\gamma^{cs}(x)) \subset W_\varepsilon^{cs}(f^n(x))$ for all $n \geq 0$.*

2. *For every $\gamma \leq \gamma(\delta_0)$, either $\ell(f^{-n}(W_\gamma^{cu}(x)) \to 0$ or $x \in W^u(p)$ for some $p \in \text{Per}(f_{/\Lambda})$ and $p \in W_\gamma^{cu}(x)$. This is the same for central stable manifolds: either $\ell(f^n(W_\gamma^{cs}(x))) \to 0$ or $x \in W^s(p)$ for some $p \in \text{Per}(f_{/\Lambda})$ and $p \in W_\gamma^{cs}(x)$.*

*In particular for every $x \in \Lambda$ there exists $\gamma_x$ such that $\ell(f^{-n}(W_{\gamma_x}^{cu}(x))) \to 0$, and it is the same for $W^{cs}$.*

*Proof.* We shall prove the lemma only for $W^{cu}$ since the other case is analogous.

1) Take some $\varepsilon \leq \delta_0$. Recall from Corollary 3.2 the existence of $\delta$ ($\delta < \varepsilon$) such that if $y \in W_\varepsilon^{cu}(x)$ and $\text{dist}(f^{-j}(x), f^{-j}(y)) \leq \delta$ for $0 \leq j \leq n$, then $f^{-j}(y) \in W_\varepsilon^{cu}(f^{-j}(x))$ for $0 \leq j \leq n$.

Assume the thesis of the lemma is false. Then there exist sequences $\gamma_n \to 0$, $x_n \in \Lambda$ and $m_n \to \infty$ such that, for $0 \leq j \leq m_n$,

$$\ell(f^{-j}(W_{\gamma_n}^{cu}(x_n))) \leq \delta$$

and

$$\ell(f^{-m_n}(W_{\gamma_n}^{cu}(x_n))) = \delta.$$

Letting $I_n = f^{-m_n}(W^{cu}_{\gamma_n}(x_n))$ we can assume (taking a subsequence if necessary) that $I_n \to I$ and $f^{-m_n}(x_n) \to z$, $z \in \Lambda$, $z \in \bar{I}$ (the closure of $I$).

Now, we have that $\ell(f^n(I)) \leq \delta$ for all positive $n$, and since $I \subset W^{cu}_\varepsilon(z)$, we conclude that $I$ is a $\delta$-$E$-interval. Thus, $\omega(z)$ is a periodic orbit $p$ because $z \in \bar{I}$. Since $z \in \Lambda$ we conclude that $p \in \Lambda$, and therefore $p$ is a hyperbolic periodic point. Hence $z \in W^s(p)$, and, at least, one component of $W^u(p) - \{p\}$ has length less than $\delta$.

In case
$$f^{-m_n}(W^{cu}_{\gamma_n}(x_n)) \cap W^s(p) \neq \emptyset$$
we get a contradiction with the inclination lemma (or $\lambda$-lemma, see [P2]) because this intersection is transversal and
$$\ell(f^{m_n}(f^{-m_n}(W^{cu}_{\gamma_n}(x_n)))) = \ell(W^{cu}_{\gamma_n}(x_n)) \to 0.$$

On the other hand, if
$$f^{-m_n}(W^{cu}_{\gamma_n}(x_n)) \cap W^s(p) = \emptyset$$
it follows, for sufficiently large $n$, that $\omega(f^{-m_n}(x_n))$ is the other endpoint of the component $W^u(p) - \{p\}$ having length less than $\delta$. By the previous lemma, it is a sink or a nonhyperbolic periodic point. This is a contradiction to the fact that $\omega(f^{-m_n}(x_n)) \subset \Lambda$, which completes the proof of 1).

2) Take $\gamma \leq \gamma(\delta_0)$ from part 1). Assume that there exists a point $x \in \Lambda$ such that $\ell(f^{-n}(W^{cu}_\gamma(x)))$ does not go to zero. By 1), we have $\ell(f^{-n}(W^{cu}_\gamma(x))) \leq \delta_0$ for all $n \geq 0$. Since we are assuming that $\ell(f^{-n}(W^{cu}_\gamma(x))$ does not converge to zero, there exist $\eta > 0$ and a sequence $n_k \to \infty$ such that
$$\ell(f^{-n_k}(W^{cu}_\gamma(x))) > \eta.$$

Letting $I_{n_k} = f^{-n_k}(W^{cu}_\gamma(x))$ we can assume that $I_{n_k} \to I$ and $f^{-n_k}(x) \to z \in \bar{I}, z \in \Lambda$. As we did in 1), we get that $I$ is a $\delta_0$-$E$-interval, and so $z \in W^s(p)$, for some periodic point $p$ in $\Lambda$.

If $z \in \text{int}(I)$, then, since $I$ is transversal to $W^s(p)$, it follows, by the inclination lemma, that $\ell(W^u(p)) \leq \delta_0$. Hence, as $f^{-n_k}(x) \to z$, we conclude that $f^{-n_k}(x) \in W^s(p)$. Otherwise, since both components of $W^u(p) - \{p\}$ have length less than $\delta_0$, we get that $\omega(x)$ is sink or a nonhyperbolic periodic point. If $z \neq p$, we may assume that the $f^{-n_k}(x)$ are contained in a fundamental domain of $W^s(p)$. Since $n_k \to \infty$ and the stable manifold does not have self intersection, we get a contradiction. Thus, $z = p$. It follows that $f^{-n_k}(x)$ belongs to the local stable manifold of $p$. If this is not the case, we get that $\omega(x)$ is a sink or a nonhyperbolic periodic point. But the only way that $f^{-n_k}(x)$ belongs to the local stable manifold of $p$ is when $x = p$. This proves 2) in case $z \in \text{int}(I)$.

Assume, now, that $z \notin \text{int}(I)$. Again, the inclination lemma implies that one component of $W^{\mathrm{u}}(p) - \{p\}$ has length less than $\delta_0$. As in 1), the case

$$f^{-n_k}(W^{cu}_\gamma(x)) \cap W^{\mathrm{s}}(p) = \emptyset$$

leads to a contradiction. So

$$f^{-n_k}(W^{cu}_\gamma(x)) \cap W^{\mathrm{s}}(p) \neq \emptyset.$$

By the inclination lemma, the facts that $\ell(f^j(f^{-n_k}(W^{cu}_\gamma(x))) \leq \delta_0$, $0 \leq j \leq n_k$, together with $f^{-n_k}(x) \to z$ imply that $x \in W^{\mathrm{u}}(p)$. Moreover, $x \in W^{\mathrm{u}}_{\delta_0}(p)$. Hence $p \in W^{cu}_\gamma(x)$, because otherwise, as $f^{-n}(x) \to p$ and $\ell(f^{-n}(W^{cu}_\gamma(x)))$ does not converge to zero, we have that $f^{-n_k}(W^{cu}_\gamma(x)) \cap W^{\mathrm{s}}(p) \neq p$ for any large $n_k$, and in particular we get $I \subset W^{\mathrm{s}}(p)$, a contradiction. This completes the proof of 2). $\square$

*Remark* 3.3. Let $\gamma_0 \leq \frac{\gamma(\delta_0)}{2}$. We may assume that there is coherence between $W^{cu}_\gamma(x)$, $x \in W^{\mathrm{u}}_\gamma(p)$ and $W^{cu}_\gamma(p)$ where $p \in \Lambda$ is a periodic point, and $\gamma \leq \gamma_0$. Furthermore, by the dynamics properties of $W^{cu}_\gamma(y)$, where $y$ does not belong to the unstable manifold of size $\gamma$ of any periodic point in $\Lambda$, it is possible to prove the coherence of the local central unstable manifolds; i.e., $W^{cu}_\gamma(x) \cap W^{cu}_\gamma(y)$ is a relative open set for both, for any points $x, y \in \Lambda$. The same applies to the central stable manifolds.

3.4. *Boxes and distortion.* Recall that for any $0 < a \leq 1$ we have an admissible neighborhood $V_a = V_a(\Lambda)$ where we have defined the central unstable and stable $a$-cone fields. From now on, we fix $V = V_1$ and $V_a$ for some $0 < a < 1$. Remember that an interval $I$ is $a$-transversal to the $E$-direction if $I \subset V_a$ and $T_xI$ lies in the central unstable $a$-cone. In the sequel we will say only transversal instead of $a$-transversal.

We shall take $\gamma < \gamma_0 \leq \frac{\gamma(\delta_0)}{2}$ ($\delta_0$ from Lemma 3.3.2) so small that, for any $x \in \Lambda$, $W^{cu}_\gamma(x)$ is transversal to the $E$-direction, and $W^{cs}_\gamma(x)$ is transversal to the $F$-direction.

Furthermore, from now on, we fix $0 < \lambda < \lambda_1 < \lambda_2 < \lambda_3 < 1$ and $c > 0$, where $\lambda$ is as in the definition of dominated splitting, $\lambda_1 = \lambda^{\frac{1}{2}}$, and $\lambda_1(1+c) < \lambda_2$, $\lambda_2(1+c) < \lambda_3$. We shall take $\delta_2 > 0$ and assume that $a$ is small enough such that

$$(1+c)^{-1} < \frac{\|Df_{/\tilde{E}(x)}\|}{\|Df_{/\tilde{E}(y)}\|} < 1+c$$

and

$$(1+c)^{-1} < \frac{\|Df^{-1}_{/\tilde{F}(x)}\|}{\|Df^{-1}_{/\tilde{F}(y)}\|} < 1+c$$

whenever $\tilde{E}$ (resp. $\tilde{F}$) lies in the central stable (resp. unstable) $a$-cone and $\text{dist}(x, y) \leq \delta_2$.

*Definition* 3. Let $x \in \Lambda$, and let $J$ be an open interval, $x \in J \subset W_\gamma^{cu}(x)$, and let $T$ be an open interval, $x \in T \subset W_\varepsilon^{cs}(x)$. A box $B_T(J)$ for $x$, is an open rectangle centered at $x$, having $J$ and $T$ as its axes, and having the boundary transversal to the $E$- and $F$-direction.

More precisely
$$B_T(J) = \text{int}(h([-1,1]^2))$$
where $h : [-1,1]^2 \to M$ is a homeomorphism, such that
$$h(\{0\} \times [-1,1]) = \bar{J}$$
and
$$h([-1,1] \times \{0\}) = \bar{T}.$$
Moreover, if we define the central stable boundary as
$$\partial^{cs}(B_T(J)) = h(\{-1,1\} \times [-1,1]),$$
we require that the two components (intervals) of it are transversal to the $E$-direction. For the central unstable boundary
$$\partial^{cu}(B_T(J)) = h([-1,1] \times \{-1,1\})$$
we require transversality to the $F$-direction. Also, for any $y \in \Lambda$, we require that $W_\gamma^{cu}(y) \cap \partial^{cs}(B_T(J))$ be a relative open set for both, and the same for $W_\gamma^{cs}(y) \cap \partial^{cu}(B_T(J))$.

If $B_T(J)$ is a box and $x \in T' \subset T$, we say that $B_{T'}(J)$ is a subbox of $B_T(J)$ if it is a box, $B_{T'}(J) \subset B_T(J)$ and $\partial^{cu}(B_{T'}) \subset \partial^{cu}(B_T(J))$.

*Remark* 3.4. These *boxes* already exist when $J$ and $\varepsilon$ are arbitrarily small, and the boundary is part of central stable and unstable manifolds of points near $x$.

*Note.* In order to simplify the notation, we will write $B_\varepsilon(J)$ instead of $B_T(J)$, $T \subset W_\varepsilon^{cs}(x)$. Also, for $\varepsilon' < \varepsilon$, we write $B_{\varepsilon'}(J)$ instead of $B_{T'}(J)$, $T' \subset T \cap W_{\varepsilon'}^{cs}(x)$.

Moreover, ordering $J$ and $W_\varepsilon^{cs}(x)$ in some way, we denote $J^+ = \{y \in J : y > x\}$, $J^- = \{y \in J : y < x\}$ and the same for $W_\varepsilon^{cs}(x)$. Also we shall denote by $B_\varepsilon(J^+)$ (say the upper part of the box) the connected component of $B_\varepsilon(J) - W_{\varepsilon(x)}^{cs}$ which contains $J^+$, and by $B_\varepsilon(J^-)$ (the bottom one) the one containing $J^-$. In an analogous way we define the left part of the box $(B_\varepsilon^-(J))$ and the right one $(B_\varepsilon^+(J))$.

Let $B_\varepsilon(J)$ be a box and $y \in B_\varepsilon(J) \cap \Lambda$. We define
$$J(y) = W_\gamma^{cu}(y) \cap B_\varepsilon(J).$$

*Remark* 3.5. Notice that, if the box is small enough, $J(y)$ is always an interval, transversal to the $E$-direction, whose endpoints are in $\partial^{cu}(B_\varepsilon(J))$.

Furthermore, as a consequence of Lemma 3.3.2, given $\delta > 0$, if the box is small enough,
$$\ell(f^{-n}(J(y))) \leq \delta$$
for any $y \in B_\varepsilon(J) \cap \Lambda$ and $n \geq 0$.

*Definition* 4. A box $B_\varepsilon(J)$ has distortion (or *cu*-distortion) $C$ if for any two intervals $J_1$, $J_2$ in $B_\varepsilon(J)$, transversal to the $E$-direction, whose endpoints are in $\partial^{cu}(B_\varepsilon(J))$ the following holds:
$$\frac{1}{C} \leq \frac{\ell(J_1)}{\ell(J_2)} \leq C.$$

*Remark* 3.6. If a box has distortion $C$, then, for any $y, z \in B_\varepsilon(J) \cap \Lambda$,
$$\frac{1}{C} \leq \frac{\ell(J(z))}{\ell(J(y))} \leq C.$$

Notice that in order to guarantee distortion $C$ on a box $B_\varepsilon(J)$, it is sufficient to find a $C^1$ foliation close to the $E$-direction in the box, such that, for any two intervals $J_1$, $J_2$ (taken as in the definition of distortion),
$$\frac{1}{C} \leq \|\Pi'\| \leq C$$
holds, where $\Pi = \Pi(J_1, J_2)$ is the projection along the foliation between these intervals.

Consider $B_\varepsilon(J)$ a box, and take a $C^1$-vector field $X$ in $B_\varepsilon(J)$, $C^0$-close to the $E$-direction ($X(x)$ lies in the central stable $b$-cone for $b$ small), and such that, for $x \in \partial^{cu}(B_\varepsilon(J))$, $X(x) \in T_x \partial^{cu}(B_\varepsilon(J))$. Consider the foliation $\mathcal{F}^{cs}$ (or the flux) generated by this vector field. For any $x \in B_\varepsilon(J)$ let $\mathcal{F}^{cs}(x)$ be the leaf passing through $x$. Notice that there exists $C$ such that
$$\frac{1}{C} \leq \|\Pi'\| \leq C$$
where $\Pi = \Pi(J_1, J_2)$ is the projection along this foliation between two intervals transversal to the $E$-direction.

The following lemma will be useful in the sequel.

LEMMA 3.4.1. *Let $B_\varepsilon(J)$, $\mathcal{F}^{cs}$, and $C$ be as above. There exist $\tau$ and $\alpha$ such that if $z \in B_\varepsilon(J) \cap \Lambda$ and for some $n > 0$ there is a box $B(n) = B_\varepsilon(f^{-n}(J(z)))$ satisfying*

1. $f^n(B(n)) \subset B_\varepsilon(J)$ *and* $f^n(\partial^{cu}(B(n)) \subset \partial^{cu}(B_\varepsilon(J))$,

2. $f^j(B(n))$ *has diameter less than $\tau$, $0 \leq j \leq n$,*

3. *there exists $K$ such that, for every $x \in B(n)$,*

$$\sum_{j=0}^{n} \ell(f^j(\mathcal{F}_n^{cs}((x))))^\alpha \leq K$$

(where $\mathcal{F}_n^{cs}(x)$ denotes the connected component of $f^{-n}(\mathcal{F}^{cs}(f^n(x))) \cap B(n)$ which contains $x$),

*then there exists $C_1 = C_1(C, K)$ such that $B(n)$ has distortion $C_1$.*

*Proof.* Consider the foliation $\mathcal{F}_n^{cs}$ in $B(n) = B(f^{-n}(J(z)))$ which is the iteration of the original one under $f^{-n}$. Let $J_1^n$, $J_2^n$ be any two intervals transversal to the $E$-direction, with endpoints in $\partial^{cu}(B(n))$ and let $\Pi_n = \Pi(J_1^n, J_2^n)$ be the projection along the foliation from $J_1^n$ to $J_2^n$. Notice that $J_1 = f^n(J_1^n)$ and $J_2 = f^n(J_2^n)$ are also two intervals in $B_\varepsilon(J)$ transversal to the $E$-direction with endpoints in $\partial^{cu}(B_\varepsilon(J))$. Let $\Pi = \Pi(J_1, J_2)$ be the projection between these two intervals along the foliation $\mathcal{F}^{cs}$. Thus we have

$$\frac{1}{C} \leq \|\Pi'\| \leq C.$$

The lemma is proved if we show that there exists $C_1 = C_1(C, K)$ such that

$$\frac{1}{C_1} \leq \|\Pi'_n\| \leq C_1.$$

For a point $x \in f^j(J_i^n)$, $i = 1, 2$ set $\tilde{F}(x) = T_x f^j(J_i^n)$, $0 \leq j \leq n$.

By the equality

$$\Pi_n \circ f_{/J_1}^{-n} = f^{-n} \circ \Pi$$

we conclude, for $y \in J_1$, that

$$\|\Pi'_n(f^{-n}(y))\|.\|Df_{/\tilde{F}(y)}^{-n}\| = \|Df_{/\tilde{F}(\Pi_1(y))}^{-n}\|.\|\Pi'(y)\|.$$

Hence

$$\|\Pi'_n(f^{-n}(y))\| = \frac{\|Df_{/\tilde{F}(\Pi(y))}^{-n}\|}{\|Df_{/\tilde{F}(y)}^{-n}\|}.\|\Pi'(y)\|.$$

Thus, to finish the proof of the lemma it suffices to find $K_1$ such that

$$\frac{1}{K_1} \leq \frac{\|Df_{/\tilde{F}(\Pi(y))}^{-n}\|}{\|Df_{/\tilde{F}(y)}^{-n}\|} \leq K_1,$$

which is the same, when $x = f^{-n}(y)$, as

$$\frac{1}{K_1} \leq \frac{\|Df_{/\tilde{F}(x)}^{n}\|}{\|Df_{/\tilde{F}(\Pi_n(x))}^{n}\|} \leq K_1.$$

With the same arguments as in [Sh, pp. 45, 46], it is possible to prove that there exist $\tau > 0$ and $\alpha > 0$ such that

$$\left| \|Df_{/\tilde{F}(f^j(w_1))}\| - \|Df_{/\tilde{F}(f^j(w_2))}\| \right| \leq \eta^j D + \operatorname{dist}(f^j(w_1), f^j(w_2))^\alpha$$

for some constants $0 < \eta < 1$ and $D$ whenever $\tilde{F}$ lies in the central unstable $a$-cone and $\operatorname{dist}(f^j(w_1), f^j(w_2)) \leq \tau$, $0 \leq j \leq n$. (This is, roughly speaking, a consequence of the fact that the distribution $F$ is $\alpha$-Hölder and any other direction converges exponentially fast to $F$.)

Therefore, it follows that

$$\frac{\|Df^n_{/\tilde{F}(x)}\|}{\|Df^n_{/\tilde{F}(\Pi_n(x))}\|} \leq \exp\left( \frac{D}{1-\eta} + \sum_{j=0}^n \operatorname{dist}(f^j(x), f^j(\Pi_n(x)))^\alpha \right).$$

Since $x$ and $\Pi_n(x)$ belong to $\mathcal{F}^{cs}_n(x)$, we conclude that

$$\sum_{j=0}^n \operatorname{dist}(f^j(x), f^j(\Pi_n(x)))^\alpha \leq \sum_{j=0}^n \ell(f^j(\mathcal{F}^{cs}_n(x)))^\alpha \leq K.$$

Thus

$$\frac{\|Df^n_{/\tilde{F}(x)}\|}{\|Df^n_{/\tilde{F}(\Pi_n(x))}\|} \leq \exp(\frac{D}{1-\eta} + K).$$

Finally, taking $K_1 = \exp(\frac{D}{1-\eta} + K)$, we have that $C_1 = CK_1$ satisfies the desired property. □

If we take the mentioned vector field in $B_\varepsilon(J)$ close enough to the $E$-direction we have the following:

COROLLARY 3.5. *There exists $\tau$ ($< \delta_2$) such that if for some $z \in B_\varepsilon(J) \cap \Lambda$ and $n > 0$ we have*

$$\|Df^j_{/E(f^{-n}(z))}\| < \lambda_1^j \text{ for } 0 \leq j \leq n$$

*and $B(f^{-n}(J(z)))$ is a box satisfying*

a) $f^n(B(n)) \subset B_\varepsilon(J)$ *and* $f^n(\partial^{cu}(B(n)) \subset \partial^{cu}(B_\varepsilon(J))$,

b) $\operatorname{diam}(f^j(B(f^{-n}(J(z))))) \leq \tau$,

*then, there exists $C_1$ such that $B(f^{-n}(J(z)))$ has distortion $C_1$.*

3.5. *More dynamical properties on $W^{cu}$ and $W^{cs}$.* In this section we shall prove that, if $\Lambda$ is transitive and every proper compact invariant subset is hyperbolic, the central unstable (stable) manifolds are dynamically defined; that is, they are unstable (stable) manifolds.

The next lemma is classical in one dimensional dynamics (see for example [dMS]) and the proof is left to the reader. However, since the diffeomorphism $f$ is of class $C^2$, so are the center unstable manifolds, and moreover these center unstable manifolds vary continuously in the $C^2$ topology. It is important to remark that there is a uniform Lipschitz constant $K_0$ of $\log(Df)$ along these manifolds.

LEMMA 3.5.1. *There exists $K_0$ such that for all $x \in \Lambda$ and $J \subset W_\gamma^{cu}(x)$, for all $z, y \in J$ and $n \geq 0$:*

1. $\dfrac{\|Df^{-n}_{/\tilde{F}(y)}\|}{\|Df^{-n}_{/\tilde{F}(z)}\|} \leq \exp(K_0 \sum_{i=0}^{n-1} \ell(f^{-i}(J)))$;

2. $\|Df^{-n}_{/\tilde{F}(x)}\| \leq \dfrac{\ell(f^{-n}(J))}{\ell(J)} \exp(K_0 \sum_{i=0}^{n-1} \ell(f^{-i}(J)))$.

The next lemma is the main result of this section. Moreover, the proof of it will be a model for the proof of the Main Lemma.

LEMMA 3.5.2. *Assume that $\Lambda$ is transitive and every proper compact invariant subset is hyperbolic. Then, either*

$$\ell(f^{-n}(W_\gamma^{cu}(x))) \to 0 \text{ as } n \to \infty \text{ for all } x \in \Lambda$$

*or $F$ is expanding (i.e., $\|Df^{-n}_{/F(x)}\| \to 0$ for all $x \in \Lambda$).*

*Proof.* Notice that, if $\Lambda$ is a periodic orbit there is nothing to prove. Thus, let us assume that $\Lambda$ is not a periodic orbit.

Now, suppose that, for some $x \in \Lambda$, $\ell(f^{-n}(W_\gamma^{cu}(x)))$ does not go to zero. Then, by Lemma 3.3.2, there exists a periodic point $p \in \Lambda$ having one component of $W^u(p) - \{p\}$ with length less than $\delta_0$. By Lemma 3.3.1, the endpoint of this component different from $p$ is a sink or a nonhyperbolic periodic point.

Since $\Lambda$ is transitive (and is not a periodic orbit), there exists $x_0 \in \Lambda$ such that $x_0 \in W^s(p) - p$. Hence, there is a small neighborhood $U$ of $x_0$, such that the points of this neighborhood at one side of $W_\varepsilon^{cs}(x_0)$ ($\subset W^s(p)$) (say the "upper" one) have their $\omega$-limit set as the other endpoint (i.e. different from $p$) of the component of $W^u(p) - \{p\}$ with length less than $\delta_0$.

Let $B_\varepsilon(J)$ be a box for $x_0$, $B_\varepsilon(J) \subset U$. Then, we may assume that

$$f^n(B_\varepsilon(J^+)) \cap B_\varepsilon(J^+) = \emptyset \text{ for all } n \geq 1$$

and so

$$f^{-n}(B_\varepsilon(J^+)) \cap B_\varepsilon(J^+) = \emptyset \text{ for all } n \geq 1.$$

We shall assume that $B_\varepsilon(J)$ is so small that, for any $y \in B_\varepsilon(J) \cap \Lambda$ and $n \geq 0$,

$$\ell(f^{-n}(J(y))) \leq \delta$$

where $\delta$ and $\varepsilon$ are taken in such a way that $\delta + \varepsilon < \tau$ ($\tau$ from Corollary 3.5).

For $z \in B_\varepsilon(J) \cap \Lambda$, let $J^+(z) = J(z) \cap B_\varepsilon(J^+)$. Then, it follows that
$$f^{-n}(J^+(z)) \cap B_\varepsilon(J^+) = \emptyset \text{ for all } n \geq 1.$$

Take $\varepsilon_1 < \frac{\varepsilon}{4}$ and consider $B_{\varepsilon_1}(J)$ a subbox. We claim that there exists $K$ such that
$$\sum_{j \geq 0} \ell(f^{-j}(J^+(z))) \leq K$$

for any $z \in B_{\varepsilon_1}(J) \cap \Lambda$. For this, consider a $C^1$-foliation in $B_\varepsilon(J)$ as we did in Section 3.4, and let $C$ be the distortion of $B_\varepsilon(J^+)$ with respect to this foliation. Let $C_1 = C_1(C)$ be as in Corollary 3.5.

Let $z \in B_{\varepsilon_1}(J)$ and consider (if it exists) the sequence of integers $0 = m_0 < m_1 < m_2 < \cdots < m_i < \cdots$ such that
$$\|Df^j_{/E(f^{-m_i}(z))}\| \leq \lambda_2^j, \ 0 \leq j \leq m_i.$$

Now, for every $i$, take a box $B(m_i) = B_{\varepsilon_1}(f^{-m_i}(J^+(z)))$ such that $f^{m_i}(B(m_i)) \subset B_\varepsilon(J)$ and $f^{m_i}(\partial^{cu}(B(m_i))) \subset \partial^{cu}(B_\varepsilon(J))$. Then, it follows, from Corollary 3.5, that $B(m_i)$ has distortion $C_1$.

Now, observe that, for $i \neq j$, $B(m_i) \cap B(m_j) = \emptyset$. Otherwise, assuming $m_j > m_i$, we get $f^{-(m_j - m_i)}(J^+(z)) \cap B_\varepsilon(J) \neq \emptyset$ which is a contradiction. Take $\varepsilon_2 < \frac{\varepsilon_1}{4}$. For any $y \in \Lambda$ consider a box $B_{\varepsilon_2}(W^{cu}_{\frac{\gamma}{2}}(y))$. The union of these boxes covers $\Lambda$ and so, since $\Lambda$ is compact, we can cover it with a finite union of such boxes, say $B(y_1), \ldots, B(y_n)$. Thus, $f^{-m_i}(z)$ belongs to one of these boxes, say $B(y_k)$ (if it belongs to more than one we choose it in an arbitrary way). Let $J(m_i) = B(m_i) \cap W^{cu}_\gamma(y_k)$. Now we have, for $i \neq j$, that $J(m_i) \cap J(m_j) = \emptyset$ and (if the adapted box $B_\varepsilon(J)$ is small enough) $J(m_i)$ is an interval in $B(m_i)$ as the definition of distortion requires.

Since $B(m_i)$ has distortion $C_1$, we conclude
$$\frac{1}{C_1} \leq \frac{\ell(J(m_i))}{\ell(f^{-m_i}(J^+(z)))} \leq C_1.$$

Then
$$\sum_i \ell(f^{-m_i}(J^+(z))) \leq C_1 \sum_i \ell(J(m_i)) \leq C_1 \sum_k \ell(W^{cu}_\gamma(y_k)) = K_1.$$

Now, we must control the sum between consecutive $m_i'$s. Let $N = N(\lambda_1, \lambda_2)$ from Lemma 3.0.2 and consider $K_2 = \sup\{\|Df^j\| : 1 \leq j \leq N\}$. There are two possibilities: $m_{i+1} - m_i < N$ or $m_{i+1} - m_i \geq N$. If $m_{i+1} - m_i < N$, then
$$\sum_{j=m_i}^{m_{i+1}-1} \ell(f^{-j}(J^+(z))) \leq NK_2 \ell(f^{-m_i}(J^+(z))).$$

On the other hand, if $m_{i+1} - m_i \geq N$, then
$$\|Df^j_{/E(f^{-m_i-j}(z))}\| \geq \lambda_1^j \text{ for } N \leq j \leq m_{i+1} - m_i.$$

Thus, by the dominated splitting,
$$\|Df^{-j}_{/F(f^{-m_i}(z))}\| \leq \lambda_1^j \text{ for } N \leq j \leq m_{i+1} - m_i.$$

Then, by Lemma 3.0.2 and Corollary 3.1, there exist $\tilde{n}_i, \tilde{n}_i - m_i < N$ such that
$$\|Df^{-j}_{/F(f^{-\tilde{n}_i}(z))}\| \leq \lambda_2^j \text{ for } 0 \leq j \leq m_{i+1} - \tilde{n}_i$$

and so, for any $y \in f^{-\tilde{n}_i}(J^+(z))$ we have, setting $\tilde{F}(y) = T_y f^{-\tilde{n}_i}(J(z))$, that
$$\|Df^{-j}_{/\tilde{F}(y)}\| \leq \lambda_3^j \text{ for } 0 \leq j \leq m_{i+1} - \tilde{n}_i.$$

Hence
$$\begin{aligned}
\sum_{j=m_i}^{m_{i+1}-1} \ell(f^{-j}(J^+(z))) &\leq \sum_{j=m_i}^{\tilde{n}_i-1} \ell(f^{-j}(J^+(z))) + \sum_{j=\tilde{n}_i}^{m_{i+1}-1} \ell(f^{-j}(J^+(z))) \\
&\leq NK_2 \ell(f^{-m_i}(J^+(z))) \\
&\quad + \sum_{j=0}^{m_{i+1}-\tilde{n}_i-1} K_2 \ell(f^{-m_i}(J^+(z))) \lambda_3^j \\
&\leq \left(NK_2 + K_2 \frac{1}{1-\lambda_3}\right) \ell(f^{-m_i}(J^+(z))).
\end{aligned}$$

Therefore
$$\begin{aligned}
\sum_{j\geq 0} \ell(f^{-j}(J^+(z))) &= \sum_i \sum_{j=m_i}^{m_{i+1}-1} \ell(f^{-j}(J^+(z))) \\
&\leq \left(NK_2 + K_2 \frac{1}{1-\lambda_3}\right) \sum_i \ell(f^{-m_i}(J^+(z))) \\
&\leq \left(NK_2 + K_2 \frac{1}{1-\lambda_3}\right) K_1 = K_3.
\end{aligned}$$

Finally, if the sequence $m_i's$ does not exist, the same argument shows that
$$\begin{aligned}
\sum_{j\geq 0} \ell(f^{-j}(J^+(z))) &\leq \left(NK_2 + K_2 \frac{1}{1-\lambda_3}\right) \ell(J^+(z)) \\
&\leq \left(NK_2 + K_2 \frac{1}{1-\lambda_3}\right) L = K_4
\end{aligned}$$

where $L = \sup\{\ell(J(z)) : z \in B_\varepsilon(J) \cap \Lambda\}$. Taking $K = \max\{K_3, K_4\}$ we conclude the proof of our claim.

Now, as in Schwartz's proof of the Denjoy theorem ([Sch]), we conclude that for all $z \in B_{\varepsilon_1}(J) \cap \Lambda$ there exists $J_1(z), J^+(z) \subset J_1(z) \subset J(z)$ such that the length of $J_1(z) - J^+(z)$ is bounded away from zero (independently of $z$) and such that, for some $\tilde{K}$,
$$\sum_{n=0}^\infty \ell(f^{-n}(J_1(z))) < \tilde{K}.$$

In particular, by Lemma 3.5.1, it follows, for any $y \in J_1(z)$, that
$$\|Df^{-n}_{/\tilde{F}(y)}\| \to 0 \text{ as } n \to \infty.$$

Indeed, for any $y, w \in J_1(z)$,
$$\frac{\|Df^{-n}_{/\tilde{F}(y)}\|}{\|Df^{-n}_{/\tilde{F}(w)}\|} \leq \exp K_0 \sum_{j=0}^{n-1} \ell(f^{-j}(J_1(z))) \leq \exp K_0 \tilde{K}.$$

Then, for $y \in J(z)$,
$$\|Df^{-n}_{/\tilde{F}(y)}\| \leq \frac{\ell(f^{-n}(J_1(z)))}{\ell(J_1(z))} \exp K_0 \tilde{K}$$

and so
$$\sum_{n=0}^{\infty} \|Df^{-n}_{/\tilde{F}(y)}\| \leq \frac{\tilde{K}}{\ell(J_1(z))} \exp K_0 \tilde{K} < \infty$$

implying that
$$\|Df^{-n}_{/\tilde{F}(y)}\| \to 0 \text{ as } n \to \infty.$$

Consider
$$B_1 = \bigcup_{z \in B_{\varepsilon_1}(J) \cap \Lambda} J_1(z).$$

Notice that $B_1$ is a neighborhood of $x_0$ in $\Lambda$.

To finish the proof of the lemma, take any $y \in \Lambda$. If $\alpha(y)$ is a proper subset of $\Lambda$, then it is hyperbolic. Thus
$$\|Df^{-n}_{/F(y)}\| \to 0.$$

On the other hand, if $\alpha(y) = \Lambda$, then, for some $m_0$, $f^{-m_0}(y) \in B_1$. Thus, $f^{-m_0}(y) \in J_1(f^{-m_0}(y))$. Therefore
$$\|Df^{-n}_{/F(f^{-m_0}(y))}\| \to 0$$

and so
$$\|Df^{-n}_{/F(y)}\| \to 0.$$

This completes the proof of the lemma. □

An analogous result also holds for central stable manifolds. Thus, there exists $\gamma_1$ such that, for all $x \in \Lambda$, as $n \to \infty$,
$$\ell(f^{-n}(W^{cu}_{\gamma_1}(x))) \to 0$$

and
$$\ell(f^n(W^{cs}_{\gamma_1}(x))) \to 0.$$

In the sequel, we shall assume $\gamma \leq \gamma_1$. Now, we have:

COROLLARY 3.6. *Given $\varepsilon > 0$, there exists $n_0$ such that, for any $x \in \Lambda$ and $n \geq n_0$,*
$$\ell(f^{-n}(W^{cu}_\gamma(x))) \leq \varepsilon$$
*and*
$$\ell(f^n(W^{cs}_\gamma(x))) \leq \varepsilon.$$

3.6. *Adapted boxes.* We shall assume that $\Lambda$ is as in the previous subsection so that the results there obtained hold.

We shall consider also the set $\tilde{\Lambda}$ as the set of points which remains in $V$ and are asymptotic to $\Lambda$; i.e.,
$$\tilde{\Lambda} = \{x \in V : f^n(x) \in V \text{ for } n \in \mathbb{Z} \text{ and } \operatorname{dist}(f^n(x), \Lambda) \to 0 \text{ as } n \to \overset{+}{-} \infty\}.$$

Notice that $\tilde{\Lambda}$ has dominated splitting since it is contained in $V$. Moreover, although $\tilde{\Lambda}$ is not compact, the arguments in Lemmas 3.3.2 and 3.5.2 apply, and so the central unstable and stable manifolds of points of $\tilde{\Lambda}$ are dynamically defined.

*Definition* 5. Let $B_\varepsilon(J)$ be a box. Recall that for $y \in B_\varepsilon(J) \cap \Lambda$ we defined $J(y) = W^{cu}_\gamma(y) \cap B_\varepsilon(J)$. We say that a box $B_\varepsilon(J)$ is $\delta$-adapted (or adapted only) if for every $y \in B_\varepsilon(J) \cap \Lambda$ the following conditions are satisfied:

1. $\ell(f^{-n}(J(y))) \leq \delta$ for all $n \geq 0$;

2. $f^{-n}(J(y)) \cap B_\varepsilon(J) = \emptyset$ or $f^{-n}(J(y)) \subset B_\varepsilon(J)$ for all $n \geq 0$.

Notice that in the preceding definition the key condition for being adapted is the second one. By Lemma 3.3.2, it is enough to take the diameter of the box sufficiently small in order to satisfy the first one.

LEMMA 3.6.1. *For every nonperiodic point $x \in \Lambda$ and for every $\delta$ arbitrarily small, there exist $\delta$-adapted boxes.*

*Proof.* Let $x$ be a nonperiodic point in $\Lambda$. As above, condition 1) is already satisfied if the diameter of the box is very small (in fact we only need the "height" of the box to be small). So, we only have to prove condition 2).

We claim that there exist $J$, $x \in J \subset W^{cu}_\gamma(x)$, such that $f^{-n}(J) \cap J = \emptyset$ for $n \geq 0$. Otherwise there exists a sequence $n_k \to \infty$ such that $f^{-n_k}(x) \in W^{cu}_\gamma(x)$, $f^{-n_k}(x) \to x$. Since $\ell(f^{-n_k}(W^{cu}_\gamma(x))) \to_k 0$ then, the $\alpha$-limit set $\alpha(x)$ is a periodic orbit, and since $f^{-n_k}(x) \to x$ we conclude that $x$ is a periodic point which is a contradiction. This completes the proof of our claim.

Therefore we can take $J$ as in the claim. We order $J$ in some way and denote $J^+ = \{y \in J : y > x\}$ and $J^- = \{y \in J : y < x\}$.

Let us prove the lemma when $x$ can be approximated by points in $\tilde{\Lambda}$ on both sides of $W_\gamma^{cs}(x)$. For that, consider the set

$$A^+ = \{y \in J^+ : \exists z \in \tilde{\Lambda}, W_\gamma^{cs}(z) \cap J^+ = \{y\}\}.$$

In a similar way we define $A^-$. Let $\varepsilon$ be any arbitrarily small number and let $B_\varepsilon(J)$ be any box, and consider a sequence $\varepsilon_k \to 0$ and a sequence of suboxes $B_{\varepsilon_{k+1}}(J) \subset B_{\varepsilon_k}(J) \subset \cdots \subset B_\varepsilon(J)$ such that $\cap_k B_{\varepsilon_k}(J) = J$. For any $k$ let $U_{y,\varepsilon_k}$ be the connected component of $W_{2\gamma}^{cs}(z) \cap B_{\varepsilon_k}(J)$ which contains $y$. Now, for $y \in A^+$, consider

$$I(y, \varepsilon_k) = \{(z, n) : z \in B_{\varepsilon_k}(J) \cap \Lambda, n > 0 \text{ and } f^{-n}(J(z)) \cap U_{y,\varepsilon_k} \neq \emptyset\}$$

We claim that for some $k$ there exists some $y \in A^+$ such that $\#I(y, \varepsilon_k) < \infty$ (and the same for $A^-$). In this case it is not difficult to see that for some $k' > k$ we have $\#I(y, \varepsilon_{k'}) = 0$, because otherwise, we could obtain, for some $n > 0$, that $f^{-n}(J) \cap J \neq \emptyset$, contradicting the way we chose $J$. Call this point $y^+$, and call $y^-$ the corresponding point for $A^-$. Then, taking $J' = (y^-, y^+) \subset J$ we construct a box $B_{\varepsilon_{k'}}(J')$ for $x$, contained in $B_\varepsilon(J)$, with $\partial^{cu}(B_{\varepsilon_{k'}}(J')) = U_{y^-,\varepsilon_{k'}} \cup U_{y^+,\varepsilon_{k'}}$. Clearly this box satisfies condition 2).

Now we may prove the claim. Arguing by contradiction, assume that $\#I(y, \varepsilon_k) = \infty$ for every $k$ and for every $y \in A^+$ or $A^-$. Without loss of generality we can assume that this holds for $A^+$.

If $\{n : \exists (z, n) \in I(y, \varepsilon_k)\}$ is bounded for some $(y, \varepsilon_k)$, then for the same $y$ and for every $k' > k$ we have the same bound. But this leads to a contradiction because for some $k'$ big enough we have either $I(y, \varepsilon_{k'}) = \emptyset$ or $f^{-n}(J) \cap J \neq \emptyset$ for some $n > 0$. Then, for every $y \in A^+$ and every $k$, the set $\{n : \exists (z, n) \in I(y, \varepsilon_k)\}$ is unbounded. In particular for every $y \in A^+$ there are sequences $z_n \in B_{\varepsilon_n}(J) \cap \Lambda$, and $m_n \to \infty$, such that

$$f^{-m_n}(J(z_n)) \cap U_{y,\varepsilon_n} \neq \emptyset.$$

Notice that $\ell(f^{-m_n}(J(z_n))) \to 0$, because otherwise, taking an accumulation point $z$ of $z_n$, we obtain that $\ell(W_\gamma^{cu}(z))$ does not go to zero, a contradiction. In particular, we have that $f^{-m_n}(z_n) \to y$, implying that $y \in \Lambda$. This means that $A^+ \subset \Lambda$. If for some $y \in A^+$ we have $y \in W^s(p)$ for some periodic point $p \in \Lambda$ we will get a contradiction. Indeed, if $y \in W^s(p)$, then $U_{y,\varepsilon_n}$ is contained in some fundamental domain of $W^s(p)$ for $n$ arbitrarily large. On the other hand, since $p \notin J$ we can take $\varepsilon_{n_0}$ sufficiently small such that the orbit of $p$ does not intersect $B_{\varepsilon_{n_0}}$. But, for $n$ large, with

$$y_n = f^{-m_n}(J(z_n)) \cap W^s(p),$$

it follows that $f^{m_n}(y_n)$ converges to the orbit of $p$ as $n \to \infty$, and so, since

$$f^{m_n}(y_n) \in f^{m_n}(f^{-m_n}(J(z_n))) = J(z_n) \subset B_{\varepsilon_n}(J) \subset B_{\varepsilon_{n_0}}(J),$$

we conclude that the orbit of $p$ intersects $B_{\varepsilon_{n_0}}(J)$, a contradiction.

Then, no $y \in A^+$ belongs to the stable manifold of a periodic point. Take some $y \in A^+$ and take $n, z_n$ and $m_n$ such that for every point $w \in J(z_n) \cap \tilde{\Lambda}$ we have $W_\gamma^{cs}(f^{-m_n}(w)) \cap J \subset J^+$ and consider the following map from $A^+$ to itself: take a point $z \in A^+$, and consider the point $w = W_\gamma^{cs}(z) \cap J(z_n)$; next take the point $f^{-m_n}(w)$; the image of the point $z$ is the point $W_\gamma^{cs}(f^{-m_n}(w)) \cap J$.

Since $\text{dist}(f^j(z), f^j(w)) \to 0$ and $\text{dist}(f^{-j}(z_n), f^{-j}(w)) \to 0$ as $j \to \infty$, the point $w$ is in $\tilde{\Lambda}$, Thus, the map is well defined.

Moreover, it is continuous and monotone as well. Now, since $A^+$ is in $J$, this implies that there exists $y_0 \in A^+$ which is a fixed point of this map. But this means that $f^{-m_n}(w_0) \in W_\gamma^{cs}(w_0)$ and since this central stable manifold is dynamically defined this implies that $\omega(w_0)$ is a periodic point and the same for $y_0$ contradicting our assumption.

Thus, we have completed the proof of the lemma when the point $x$ is accumulated by points on both sides of $W_\gamma^{cs}(x)$.

Assume now, that the point $x$ is accumulated by points in $\tilde{\Lambda}$ only on one side of $W_\gamma^{cs}(x)$ (we shall refer to them as boundary points of $\tilde{\Lambda}$). Take a small neighborhood $U(x)$ such that in this neighborhood there are points of $\tilde{\Lambda}$ only on one side of $W_\gamma^{cs}(x)$ (say the lower side), and take a small box $B_\varepsilon(J)$ contained in this neighborhood; thus $B_\varepsilon(J^+) \cap \tilde{\Lambda} = \emptyset$, where $B_\varepsilon(J^+)$ is the upper part of the box. Moreover, by the previous arguments, we can assume that this box is adapted at the bottom, meaning that, for every $y \in B_\varepsilon(J) \cap \Lambda$,

$$f^{-n}(J(y)) \cap \partial^{cu,-}(B_\varepsilon(J)) = \emptyset, \ n \geq 0$$

where $\partial^{cu,-}(B_\varepsilon(J))$ is the central unstable boundary contained in the bottom part of the box.

We shall prove then, that for any arbitrarily small $\varepsilon' < \varepsilon$, the subbox $B_{\varepsilon'}(J)$ is adapted. Arguing by contradiction, suppose that there exist a sequence $\varepsilon_n \to 0$ and subboxes $B_{\varepsilon_n}(J)$, $\cap_n B_{\varepsilon_n}(J) = J$ which are not adapted. Thus, we have sequences $y_n \in B_{\varepsilon_n}(J) \cap \Lambda$ and $m_n$ such that

$$f^{-m_n}(J(y_n)) \cap \partial^{cu,+}(B_{\varepsilon_n}(J)) \neq \emptyset.$$

The sequence $m_n$ must be unbounded; otherwise for some $n > 0$ we have $f^{-n}(J) \cap J \neq \emptyset$ which is contradiction.

But then, $\ell(f^{-m_n}(J(y_n)))$ is arbitrarily small when $n$ is large enough (remember that the central manifolds are dynamically defined). Hence, for $n$ large enough, $f^{-m_n}(J(y_n))$ is contained in the upper side of the neighborhood. Thus, $f^{-m_n}(y_n)$ is also, and it is in $\Lambda$. This is a contradiction to the way we chose $U(x)$, proving that for some $\varepsilon'$ small, the box $B_{\varepsilon'}(J)$ is adapted. $\square$

*Remark* 3.7. As a consequence of the proof of the preceding lemma, for any nonperiodic point $x$, small adapted boxes can be taken in such a way that the $\partial^{cu}(B_\varepsilon(J))$ is contained in some central stable manifold, unless the point $x$ is a boundary point of $\tilde{\Lambda}$.

Notice that, if $B_\varepsilon(J)$ is an adapted box, then any subbox is adapted also.

LEMMA 3.6.2. *Let $x \in \Lambda$ be a boundary point of $\tilde{\Lambda}$. Then, for any small adapted box $B_\varepsilon(J)$ for $x$, the side of the box which contains point of $\Lambda$ (say $B_\varepsilon(J^-)$) is adapted; i.e., for $y \in B_\varepsilon(J) \cap \Lambda$ and $n \geq 0$, either*

$$f^{-n}(J^-(y)) \subset B_\varepsilon(J^-)$$

*or*

$$f^{-n}(J^-(y)) \cap B_\varepsilon(J^-) = \emptyset$$

*where $J^-(y) = J(y) \cap B_\varepsilon(J^-)$.*

*Proof.* Let $B_\varepsilon(J)$ be a small adapted box for $x$ such that $B_\varepsilon(J^+) \cap \tilde{\Lambda} = \emptyset$. We can assume also that $\partial^{cu,-}(B_\varepsilon(J))$ is contained in some central stable manifold of a point of $\tilde{\Lambda}$ (see the proof of the previous lemma). Now, for any $y \in B_\varepsilon(J) \cap \Lambda$, the endpoints of $J^-(y)$ are in $\tilde{\Lambda}$. Thus, if $B_\varepsilon(J^-)$ is not adapted, then, for some $y \in B_\varepsilon(J) \cap \Lambda$ and $n > 0$, $f^{-n}(J^-(y)) \cap W^{cs}_\varepsilon(x) \neq \emptyset$. Since $B_\varepsilon(J)$ is adapted, $f^{-n}(J(y))$ (and so $f^{-n}(J^-(y))$) is contained in $B_\varepsilon(J)$. Therefore, one of the endpoints of $f^{-n}(J^-(y))$ is in $B_\varepsilon(J^+)$. But this point is in $\tilde{\Lambda}$ also, a contradiction. $\square$

*Remark 3.8.* When $x$ is a boundary point of $\tilde{\Lambda}$, a box adapted for it will be as in the preceding lemma, i.e., $B_\varepsilon(J^-)$. Nevertheless, abusing the notation, we shall refer to it as $B_\varepsilon(J)$. Notice that this box can be taken with the central unstable boundary contained in central stable manifolds.

*Definition 6.* Let $B_\varepsilon(J)$ be an adapted box. A map $\psi : S \to B_\varepsilon(J)$, where $S \subset B_\varepsilon(J)$, is called a return of $B_\varepsilon(J)$ associated to $\Lambda$ if:

- $S \cap \Lambda \neq \emptyset$.

- There exists $k > 0$ such that $\psi = f^{-k}_{/S}$.

- $\psi(S) = f^{-k}(S)$ is a connected component of $f^{-k}(B_\varepsilon(J)) \cap B_\varepsilon(J)$.

- $f^{-i}(S) \cap B_\varepsilon(J) = \emptyset$ for $1 \leq i < k$.

We will denote the set of returns of $B_\varepsilon(J)$ associated to $\Lambda$ by $\mathcal{R}(B_\varepsilon(J), \Lambda)$. Moreover, a return $\psi \in \mathcal{R}(B_\varepsilon(J), \Lambda)$ has $|\psi'| < \xi < 1$ if $\|Df^{-k}_{/\tilde{F}(y)}\| < \xi$ for all $y \in J(z)$, $z \in \mathrm{dom}(\psi) \cap \Lambda$, where $\psi = f^{-k}_{/\mathrm{dom}(\psi)}$ and where $\tilde{F}(y) = T_y J(z) = T_y W^{cu}_\gamma(z)$.

*Remark 3.9.* 1. If $\psi \in \mathcal{R}(B_\varepsilon(J), \Lambda)$ then $S = \mathrm{dom}(\psi)$ is a vertical strip in $B_\varepsilon(J)$ (i.e., a box such that $\partial^{cu}(S) \subset \partial^{cu}(B_\varepsilon(J))$). Notice that for every $z \in B_\varepsilon(J) \cap \Lambda$, either $J(z) \subset S$ or $J(z) \cap S = \emptyset$ holds.

2. If $\psi : S \to B_\varepsilon(J)$ is a return then $\psi(S)$ is an adapted box. Moreover, if $S \cap \partial^{cs}(B_\varepsilon(J)) = \emptyset$ then, $\psi(S)$ is a horizontal strip, i.e., a box such that $\partial^{cs}(\psi(S)) \subset \partial^{cs}(B_\varepsilon(J))$.

3. If $\psi_1 : S_1 \to B_\varepsilon(J)$ and $\psi_2 : S_2 \to B_\varepsilon(J)$ are two different returns then $S_1 \cap S_2 = \emptyset$ and $\psi(S_1) \cap \psi(S_2) = \emptyset$.

*Definition* 7. Let $B_\varepsilon(J)$ be an adapted box. We say that $B_\varepsilon(J)$ is well adapted if there exist a subbox $B_{\varepsilon'}(J)$ and two disjoint vertical strips $S_1, S_2$, such that
$$B_\varepsilon(J) - B_{\varepsilon'}(J) = S_1 \cup S_2$$
where $S_1, S_2$ satisfy either

a) $S_i \cap \Lambda = \emptyset$

or

b) $S_i$ is a domain of a return $\psi_i \in \mathcal{R}(B_\varepsilon(J), \Lambda)$, and $\psi_i(S_i)$ is a horizontal strip.

LEMMA 3.6.3. *There exist well-adapted boxes of arbitrarily small diameter.*

*Proof.* Let $B_\varepsilon(J)$ be an arbitrarily small adapted box. Observe that for any $\varepsilon' < \varepsilon$ we have that the subbox $B_{\varepsilon'}(J)$ is also an adapted box. Thus, we only have to prove that for some $\varepsilon' < \varepsilon$, $B_{\varepsilon'}(J)$ is well-adapted. We shall divide the proof in two cases:

*Case* 1.　Assume that for any $\varepsilon' \leq \varepsilon$ the subbox $B_{\varepsilon'}(J)$ satisfies:
$$B_{\varepsilon'}(J) \neq \bigcup \{J(z) : z \in \mathrm{Cl}(B_{\varepsilon'}(J)) \cap \Lambda\}.$$
In particular, this is true either for the right part of the box or for the left one. Without loss of generality, we may suppose that this is true for the right; i.e., for any $\varepsilon' \leq \varepsilon$,
$$B_{\varepsilon'}^+(J) \neq \bigcup \{J(z) : z \in \mathrm{Cl}(B_\varepsilon^+(J)) \cap \Lambda\}.$$
In particular, this holds for $\varepsilon' = \varepsilon$. Then, there exists some $w \in W_\varepsilon^{cs,+}(x) = \{y \in W_\varepsilon^{cs}(x) : y > x\}$ such that $w \notin J(z)$ for every $z \in B_\varepsilon^+(J) \cap \Lambda$. If for every $z \in \mathrm{Cl}(B_\varepsilon^+(J)) \cap \Lambda$ we have that the point $w_z = J(z) \cap W_\varepsilon^{cs,+}(x)$ is at the left of $w$. Then there exists a vertical strip $S = S_2$ at the right of $B_\varepsilon(J)$ such that $S \cap \Lambda = \emptyset$. On the other hand, if there exists $z \in \mathrm{Cl}(B_\varepsilon^+(J)) \cap \Lambda$ such that the point $w_z$ is at the right of $w$, then we consider $w_0 = \inf\{w_z : w_z > w\}$. It follows that $w_0 > w$ and then it is easy to construct a vertical strip $S \subset B_\varepsilon^+(J)$ such that $S \cap \Lambda = \emptyset$. Then it is possible to reduce the box $B_\varepsilon(J)$ having $S$ at the right boundary (for simplicity, call this new box $B_\varepsilon(J)$).

We claim now that the same holds for the left part, i.e.

$$B_\varepsilon^-(J) \neq \bigcup \{J(z) : z \in \mathrm{Cl}(B_\varepsilon^-(J)) \cap \Lambda\}.$$

Suppose this is not the case. Then

$$B_\varepsilon^-(J) = \bigcup \{J(z) : z \in \mathrm{Cl}(B_\varepsilon^-(J)) \cap \Lambda\}.$$

If there is only one return $\psi \in \mathcal{R}(B_\varepsilon(J), \Lambda)$ such that $\mathrm{dom}(\psi) \cap B_\varepsilon^-(J) \neq \emptyset$ (the set of such returns is nonempty by the transitivity of $\Lambda$), it follows that $B_\varepsilon^-(J) \subset \mathrm{dom}(\psi)$. Moreover, $\psi(B_\varepsilon^-(J)) \subset B_\varepsilon^-(J)$, because otherwise, it contradicts our assumption on the right part of the box. But then, since we may assume that $f^{-n}(J) \cap J = \emptyset$ for $n > 0$, $\mathrm{Cl}(\psi(B_\varepsilon^-(J))) \cap J = \emptyset$. Hence, no point $y \in B_\varepsilon^-(J) \cap \Lambda$ ever passes through $B_\varepsilon^-(J) - \mathrm{Cl}(\psi(B_\varepsilon^-(J)))$. This contradicts the transitivity of $\Lambda$. Therefore, there is more than one return $\psi$ such that $\mathrm{dom}(\psi) \cap B_\varepsilon^-(J) \neq \emptyset$. We see easily, applying the same arguments, that it should be more than two. Thus, there exists a return $\psi \in \mathcal{R}(B_\varepsilon(J))$, such that $\mathrm{dom}(\psi) \subset B_\varepsilon^-(J)$ and $\mathrm{Image}(\psi) = B_\psi$ is a horizontal strip. Then,

$$B_\psi = \bigcup \{J_\psi(z) : z \in \mathrm{Cl}(B_\psi) \cap \Lambda\}.$$

Since $B_\psi$ is a full horizontal strip, we also conclude that

$$B_\varepsilon(J) = \bigcup \{J(z) : z \in \mathrm{Cl}(B_\varepsilon(J)) \cap \Lambda\},$$

and so the same holds for the right part of $B_\varepsilon^+(J)$, which is a contradiction, proving our claim.

Now, we could repeat the argument used for the left part of the box as we did for the right one, and after reducing the box, we could also construct a well-adapted box.

*Case* 2. Assume now that

$$B_\varepsilon(J) = \bigcup \{J(z) : z \in \mathrm{Cl}(B_\varepsilon(J)) \cap \Lambda\}.$$

We may take (reducing the box if necessary) a return $\psi_1$ such that $S_1 = \mathrm{dom}(\psi_1) \cap B_\varepsilon^-(J) \neq \emptyset$, $\psi_1(S_1)$ is a horizontal strip and a return $\psi_2$ such that $S_2 = \mathrm{dom}(\psi_2) \cap B_\varepsilon^+(J) \neq \emptyset$ and $\psi_2(S_2)$ is a horizontal strip (we can take these returns to be different). We say that the $\psi_i$ ($i = 1, 2$) preserves or reverses the orientation when the map

$$P_i : S_i \cap W_\varepsilon^{cs}(x) \to W_\varepsilon^{cs}(x),$$

defined by $P(w) = J(\psi_i(z_w)) \cap W_\varepsilon^{cs}(x)$ where $w \in J(z_w)$, preserves or reverses the orientation of $W_\varepsilon^{cs}(x)$.

Assume that $P_1$ and $P_2$ preserve the orientation. Then, take $w_1$ a fixed point under $P_1$ and $w_2$ a fixed point under $P_2$. Then, cutting the box at $J(z_{w_1})$ and at $J(z_{w_2})$, the new box is well-adapted as we wished.

Suppose now that one preserves the orientation (say $P_1$) and the other does not. Then, take $w_1$ a fixed point under $P_1$ and $w_2$ satisfying $w_1 = P_2(w_2)$. Again, cutting the box at $J(z_{w_1})$ and at $J(z_{w_2})$, we construct a well-adapted box.

Finally, if both reverse the orientation, take $w_1$ and $w_2$ such that $P_1(w_1) = w_2$ and $P_2(w_2) = w_1$. With the same procedure as before, we conclude the proof of the lemma. □

3.7. *Proof of the Main Lemma.* In this section we shall assume that $\Lambda_0 = \Lambda$ as in the statement of the Main Lemma; i.e., $\Lambda$ is a nontrivial transitive set such that every proper compact invariant subset is hyperbolic and it is not a periodic simple curve normally hyperbolic conjugated to an irrational rotation.

We shall prove that for every $x \in \Lambda$, $\|Df^{-n}_{/F(x)}\| \to 0$ and $\|Df^n_{/E(x)}\| \to 0$ as $n \to \infty$, but we shall only show $\|Df^{-n}_{/F(x)}\| \to 0$, the other being similar when $f$ is replaced by $f^{-1}$.

LEMMA 3.7.1. *Let $B_\varepsilon(J)$ be an adapted box. There exists $K_1 = K_1(B_\varepsilon(J))$ such that if $z \in B_\varepsilon(J) \cap \Lambda$ and $f^{-i}(z) \notin B_\varepsilon(J), 1 \le i \le n$ then*

$$\sum_{i=0}^{n} \ell(f^{-i}(J(z))) \le K_1.$$

*Proof.* Let $\Lambda_1$ be the maximal invariant subset of $\Lambda$ outside of $B_\varepsilon(J)$; i.e.,

$$\Lambda_1 = \bigcap_{n \in Z} f^n(\Lambda - B_\varepsilon(J)).$$

If $\Lambda_1 = \emptyset$, then there exists $N$ such that for every $z \in B_\varepsilon(J) \cap \Lambda$ there exists $m, 1 \le m \le N$ such that $f^{-m}(z) \in B_\varepsilon(J)$. This implies (for $z$ taken according to a hypothesis of the lemma) that

$$\sum_{i=0}^{n} \ell(f^{-i}(J(z))) \le N\mathrm{diam}(M)$$

and so it is enough to take $K_1 = N\mathrm{diam}(M)$.

Assume now that $\Lambda_1 \ne \emptyset$. Since $\Lambda_1$ is a proper compact invariant subset of $\Lambda$, it is hyperbolic. Then there exist $0 < \sigma < 1$ and $n_1$ such that for $y \in \Lambda_1$ we have $\|Df^{-n}_{/F(y)}\| \le \sigma^n$, for all $n \ge n_1$. There exist $\sigma_1$, $\sigma < \sigma_1 < 1$, and $U$ a small neighborhood of $\Lambda_1$ such that if $y \in U \cap \Lambda$ and $f^{-j}(y) \in U$ for $0 \le j \le n$, where $n \ge n_1$ then

$$\|Df^{-j}_{/F(y)}\| < \sigma_1^j, \quad n_1 \le j \le n.$$

Consider now $\sigma < \sigma_1 < \sigma_2 < \sigma_3 < 1$ and $\delta_1$ such that $(1 + \delta_1)\sigma_2 < \sigma_3$ and take $N_1 \ge \max\{N(\sigma_1, \sigma_2), n_1\}$ where $N(\sigma_1, \sigma_2)$ is as in Lemma 3.0.2. We conclude that if $y \in U \cap \Lambda$ and $f^{-j}(y) \in U, 0 \leq j \leq n, n > N_1$ then there exists $m_0 < N_1$ such that

$$\|Df^{-j}_{/F(f^{-m_0}(y))}\| \leq \sigma_2^j, \ 0 \leq j \leq n - m_0.$$

Also take $\delta_2$ with the property $\frac{\|Df^{-1}_{/F(x)}\|}{\|Df^{-1}_{/F(y)}\|} \leq 1 + \delta_1$ if $\text{dist}(x,y) < \delta_2$, and $\varepsilon_1$ such that $f^{-n}(W^{cu}_{\varepsilon_1}(x)) \subset W^{cu}_{\delta_2}(f^{-n}(x))$ for all $n \geq 0$ (see Lemma 3.3.2).

Considering all these facts together, we have, for a point $y$ with the property mentioned above, that

$$\|Df^{-j}_{/F(z)}\| \leq \sigma_3^j, \ \text{for all } z \in W^{cu}_{\varepsilon_1}(f^{-m_0}(y)), \ 0 \leq j \leq n - m_0.$$

And so

$$\sum_{j=0}^{n-m_0} \ell(f^{-j}(W^{cu}_{\varepsilon_1}(f^{-m_0}(y)))) \leq \sum_{j=0}^{n-m_0} \sigma_3^j \ell(W^{cu}_{\varepsilon_1}(f^{-m_0}(y))) \leq \frac{1}{1-\sigma_3} \text{diam}(M).$$

Since $\Lambda_1$ is the maximal invariant subset of $\Lambda$ in the complement of $B_\varepsilon(J)$, there exists $m_1$ such that for every $n > m_1$,

$$\bigcap_{|j| \leq n} f^j(\Lambda - B_\varepsilon(J)) \subset U.$$

We can assume that $m_1 > N_1$. By Corollary 3.6 there exists $m_2$ such that if $n \geq m_2$ then, $f^{-n}(J(z)) \subset W^{cu}_{\varepsilon_1}(f^{-n}(z))$. Take $N = \max\{m_2 + 2N_1, 2m_1\}$ and define $K_1 = (N + \frac{1}{1-\sigma_3})\text{diam}(M)$. We will show that this $K_1$ satisfies the conclusion of the lemma.

For this, let $z \in B_\varepsilon(J) \cap \Lambda$ and assume that $f^{-j}(z) \notin B_\varepsilon(J)$ for $1 \leq j \leq n$. There are two possibilities:

- If $n \leq N$, then $\sum_{j=0}^n \ell(f^{-j}(J(z))) \leq N\delta \leq K_1$.

- If $n > N$, then there exists $n(z) < m_2 + N_1$ such that

$$f^{-n(z)}(J(z)) \subset W^{cu}_{\varepsilon_1}(f^{-n(z)}(z))$$

and

$$\|Df^{-j}_{/F(f^{-n(z)}(z))}\| \leq \sigma_2^j, \quad 0 \leq j \leq n - n(z).$$

Finally

$$\sum_{j=0}^n \ell(f^{-j}(J(z))) = \sum_{j=0}^{n(z)} \ell(f^{-j}(J(z))) + \sum_{j=0}^{n-n(z)} \ell(f^{-j}(f^{-n(z)}(J(z))))$$

$$\leq N\text{diam}(M) + \text{diam}(M)\frac{1}{1-\sigma_3} \leq K_1$$

and the proof of the lemma is finished. $\square$

Now, as a consequence of the proof of the preceding lemma, we have the following:

COROLLARY 3.7. *Let $B_\varepsilon(J)$ be an adapted box and $0 < \alpha \leq 1$. Then there exists $K_1 = K_1(B_\varepsilon(J), \alpha)$ such that if $z \in B_\varepsilon(J) \cap \Lambda$ and $f^{-i}(z) \notin B_\varepsilon(J)$, $1 \leq i \leq n$ then*

$$\sum_{i=0}^{n} \ell(f^{-i}(J(z)))^\alpha \leq K_1.$$

*Remark* 3.10. In the proof of the last lemma we did not use the fact that the box is adapted. In particular we have similar estimates for the central stable manifolds. This means that there exists $\tilde{K}_1$ such that if $z \in B_\varepsilon(J) \cap \Lambda$ and $f^j(z) \notin B_\varepsilon(J)$ for $1 \leq j \leq n$ then

$$\sum_{j=0}^{n} \ell(f^j(T(z)))^\alpha \leq \tilde{K}_1,$$

where $T(z) = W_\gamma^{cs}(z) \cap B_\varepsilon(J)$.

Furthermore, if we take a $C^1$ vector field in $B_\varepsilon(J)$ as we did in Section 3.4 close enough to the $E$-direction, then it is possible to prove, using the same arguments as in the preceding lemma, that if we take a return $\psi$ and push the foliation generated by this vector field under $\psi$ and call it $\mathcal{F}_\psi$,

$$\sum_{j=0}^{k} \ell(f^j(\mathcal{F}_\psi(x)))^\alpha \leq \tilde{K}_2$$

for some $\tilde{K}_2$ where $\psi = f^{-k}$.

In particular, taking into account Lemma 3.4.1, the following is true.

COROLLARY 3.8. *Let $B_\varepsilon(J)$ be an adapted box. Then, there exists $C_1$ such that for every return $\psi \in \mathcal{R}(B_\varepsilon(J), \Lambda)$, the adapted box Image($\psi$) has distortion $C_1$.*

LEMMA 3.7.2. *Let $B_\varepsilon(J)$ be an adapted box and assume that for every return $\psi \in \mathcal{R}(B_\varepsilon(J), \Lambda)$, $|\psi'| < \xi < 1$ for some $\xi$. Then for all $y \in B_\varepsilon(J) \cap \Lambda$ the following holds*:

$$\sum_{n \geq 0} \ell(f^{-n}(J(y))) < \infty.$$

*This implies that*

$$\|Df^{-n}_{/F(y)}\| \to_{n \to \infty} 0.$$

*Proof.* Let $y$ be as in the hypothesis of the lemma, and let $0 < n_1 < n_2 < \cdots < n_i < \cdots$ be the set of $n'$s such that $f^{-n_i}(y) \in B_\varepsilon(J)$. For every $i$, there exists $\psi_i \in \mathcal{R}(B_\varepsilon(J), \Lambda)$ such that $\psi_i = f_{/S_i}^{n_i - n_{i-1}}$ where $S_i = \text{dom}(\psi_i)$, and $f^{-n_{i-1}}(y) \in S_i$. Put $m_i = n_i - n_{i-1}$.

Let $K_1$ be as in the previous lemma. Then, for every $i$,

$$\sum_{n=0}^{m_i-1} \ell(f^{-n}(J(f^{-n_{i-1}}(y)))) < K_1 \ .$$

By Lemma 3.5.1, for $0 \le n \le m_i$ and for all $z, w \in J(f^{-n_{i-1}}(y))$ we get that

$$\frac{\|Df_{/\tilde{F}(z)}^{-n}\|}{\|Df_{/\tilde{F}(w)}^{-n}\|} \le \exp(K_0 K_1).$$

Define $K_2 = \exp(K_0 K_1)$. For every interval $L \subset J(f^{-n_{i-1}}(y))$

$$\frac{\ell(f^{-n}(L))}{\ell(L)} \le K_2 \frac{\ell(f^{-n}(J(f^{-n_{i-1}}(y))))}{\ell(J(f^{-n_{i-1}}(y)))}$$

holds, and then

$$\sum_{n=0}^{m_i-1} \ell(f^{-n}(L)) \le \frac{\ell(L) K_2}{\ell(J(f^{-n_{i-1}}(y)))} \sum_{n=0}^{m_i-1} \ell(f^{-n}(J(f^{-n_{i-1}}(y))))$$

$$\le \frac{\ell(L) K_2}{\ell(J(f^{-n_{i-1}}(y)))} K_1.$$

Setting $K_3 = K_2 K_1$ and $L = f^{-n_{i-1}}(J(y))$ we conclude

$$\sum_{n=0}^{m_i-1} \ell(f^{-n}(f^{-n_{i-1}}(J(y)))) \le K_3 \frac{\ell(f^{-n_{i-1}}(J(y)))}{\ell(J(f^{-n_{i-1}}(y)))}.$$

Since $f^{-n_{i-1}} = \psi_{i-1} \circ \cdots \circ \psi_1$, $\|Df_{/\tilde{F}(z)}^{-n_{i-1}}\| \le \xi^{i-1}$ for all $z \in J(y)$. This implies $\ell(f^{-n_{i-1}}(J(y))) \le \xi^{i-1} \ell(J(y))$ and also

$$\sum_{n=0}^{m_i-1} \ell(f^{-n}(f^{-n_{i-1}}(J(y)))) \le K_3 \xi^{i-1} \frac{\ell(J(y))}{\ell(J(f^{-n_{i-1}}(y)))}.$$

Moreover, there exists $C = C(B_\varepsilon(J))$ such that $\frac{\ell(J(z))}{\ell(J(w))} \le C$ for all $z, w \in B_\varepsilon(J) \cap \Lambda$. We conclude then

$$\sum_{n=0}^{m_i-1} \ell(f^{-n}(f^{-n_{i-1}}(J(y)))) \le K_3 C \xi^{i-1}.$$

Finally, setting $K = \max\{K_1, K_3 C\}$ we get

$$\sum_{n \ge 0} \ell(f^{-n}(J(y))) = \sum_{i=1}^{\infty} \sum_{n=0}^{m_i-1} \ell(f^{-n}(f^{-n_{i-1}}(J(y)))) \le \sum_{i=1}^{\infty} K \xi^{i-1} = \tilde{K} < \infty.$$

If the sequence of $n'_i s$ is finite or empty we can also conclude that

$$\sum_{n\geq 0} \ell(f^{-n}(J(y))) \leq \tilde{K} + K_1 < \infty.$$

In particular, as in Lemma 3.5.2, for all $z, w \in J(y)$,

$$\frac{\|Df^{-n}_{/\tilde{F}(z)}\|}{\|Df^{-n}_{/\tilde{F}(w)}\|} \leq \exp K_0 \sum_{j=0}^{n-1} \ell(f^{-j}(J(y))) \leq \exp K_0 \tilde{K}.$$

Then for $z \in J(y)$,

$$\|Df^{-n}_{/\tilde{F}(z)}\| \leq \frac{\ell(f^{-n}(J(y)))}{\ell(J(y))} \exp K_0 \tilde{K}$$

and so

$$\sum_{n=0}^{\infty} \|Df^{-n}_{/\tilde{F}(z)}\| \leq \frac{\tilde{K}}{\ell(J(y))} \exp K_0 \tilde{K} < \infty.$$

Now, taking $y = z$ in the above inequality we have

$$\sum_{n=0}^{\infty} \|Df^{-n}_{/F(y)}\| < \infty,$$

implying that

$$\|Df^{-n}_{/F(y)}\| \to 0 \text{ as } n \to \infty,$$

concluding the proof of the lemma. □

In Lemma 3.7.1 we showed that there exists a bound of the sum of the length of the iterates of the central unstable intervals of an adapted box until they return to the box. However, this bound depends on the adapted box. We will show now that we can obtain a uniform bound for adapted boxes which are images of returns of the original one.

LEMMA 3.7.3. *Let $B_\varepsilon(J)$ be a well-adapted box. There exists $K = K(B_\varepsilon(J))$ with the following property: for every $\psi \in \mathcal{R}(B_\varepsilon(J), \Lambda)$ and $z \in B_\psi = \text{Image}(\psi)$ (denoting $J_\psi(z) = J(z) \cap B_\psi$)*

$$\sum_{j=0}^{n} \ell(f^{-j}(J_\psi(z))) \leq K$$

*whenever $f^{-j}(z) \notin B_\psi, 1 \leq j \leq n$.*

*Proof.* Since $B_\varepsilon(J)$ is a well adapted box, there exist a subbox $B_{\varepsilon'}(J)$ and two disjoint vertical strip $S_1, S_2$ such that $B_\varepsilon(J) - B_{\varepsilon'}(J) = S_1 \cup S_2$ and $S_i$ is either a domain of a return or $S_i \cap \Lambda = \emptyset$.

Notice that there exists $\varepsilon_1$ such that if $y \in B_{\varepsilon'}(J) \cap \Lambda$ (i.e., $y \notin S_1 \cup S_2$), then $W^{cs}_{\varepsilon_1}(y) \subset B_\varepsilon(J)$.

Let $\psi \in \mathcal{R}(B_\varepsilon(J), \Lambda)$ and $z \in B_\psi = \text{Image}(\psi)$ such that $f^{-j}(z) \notin B_\psi, 1 \leq j \leq n$. Let $C_1$ be as in Corollary 3.8 and $C_2$ as in Corollary 3.5 (corresponding to $C_2 = C_1$).

Let $0 < n_1 < n_2 < \cdots < n_k \leq n$ be the set $\{0 < j \leq n : f^{-j}(z) \in B_\varepsilon(J)\}$. For every $n_i$ we have associated a return $\psi_i \in \mathcal{R}(B_\varepsilon(J), \Lambda)$ such that $f^{-n_i}(z) \in B_{\psi_i}$, i.e., $f^{-n_i}(z) = \psi_i(f^{-n_{i-1}}(z))$. Define $B(n_i)$ as the horizontal strip in $B_{\psi_i}$ determined by $f^{-n_i}(J_\psi(z))$, that is, the connected component of $f^{-n_i}(B_\psi) \cap B_{\psi_i}$ which contains $f^{-n_i}(z)$.

It follows that, for $i \neq j$, we have $B(n_i) \cap B(n_j) = \emptyset$. Otherwise, if for some $j > i$ we have $B(n_i) \cap B(n_j) \neq \emptyset$, then there exist $n_k$ such that $n_k = n_j - n_i$ and also $B(n_k) \cap B_\psi \neq \emptyset$. Since $B_\psi$ is adapted we conclude that $f^{-n_k}(z) \in B_\psi$ which is a contradiction because $n_k < n$.

As in the proof of Lemma 3.5.2, consider (if it exists) the sequence $0 = m_0 < m_1 < m_2 < \cdots < m_l \leq n$ such that
$$\|Df^j_{/E(f^{m_i}(z))}\| < \lambda_2^j, \ 0 \leq j \leq m_i, \ \text{for all } i = 1, \ldots, l.$$
We claim that there exists $C_4 = C_4(B_\varepsilon(J))$ such that
$$\sum_{i=0}^{l} \ell(f^{-m_i}(J_\psi(z))) \leq C_4.$$

To prove the claim, assume first that $z \notin S_1 \cup S_2$. Set $\varepsilon_2 = \frac{\varepsilon_1}{2}$, $\varepsilon_3 = \frac{\varepsilon_2}{2}$. For any point $y \in \Lambda$ consider a box (not necessarily adapted) $B_{\varepsilon_3}(W_\gamma^{cu}(y))$. Since $\Lambda$ is compact we can cover $\Lambda$ by a finite number of such boxes. We will denote these by $B_1, \ldots, B_r$. Set $C_3 = \sum_{k=1}^{r} \ell(W_{2\gamma}^{cu}(y_k))$.

For $1 \leq i \leq l$,
$$f^{m_i}(B_{\varepsilon_2}(f^{-m_i}(J_\psi(z)))) \subset B_\psi.$$

Thus, for $i \neq j$,
$$B_{\varepsilon_2}(f^{-m_i}(J_\psi(z))) \cap B_{\varepsilon_2}(f^{-m_j}(J_\psi(z))) = \emptyset.$$
Otherwise, $f^{-(m_j - m_i)}(z) \in B_\psi$ which is a contradiction since $m_j - m_i \leq n$, or which would contradict the fact that $B_\psi$ is an adapted box as well. Since $B_1, \ldots, B_r$ covers $\Lambda$, the $f^{-m_i}(z)$ belong to some of these boxes, say $B_k$ (if it belongs to more than one we choose one in an arbitrary way). Let $J_{m_i} = B_{\varepsilon_2}(f^{-m_i}(J_\psi(z))) \cap W_{2\gamma}^{cu}(y_k)$. From Corollary 3.5, for every $i$,
$$\frac{1}{C_2} \leq \frac{\ell(f^{-m_i}(J_\psi(z)))}{\ell(J_{m_i})} \leq C_2.$$
Moreover, by the fact that
$$f^{m_i}(B_{\varepsilon_2}(f^{-m_i}(J_\psi(z)))) \subset B_\psi,$$
we conclude
$$J_{m_i} \cap J_{m_j} = \emptyset.$$

Hence
$$\sum_{i=0}^{l} \ell(f^{-m_i}(J_\psi(z))) \leq \sum_{i=0}^{l} C_2 \ell(J_{m_i}) \leq C_2 C_3.$$

Assume now that $z \in S_1 \cup S_2$. Let $i_0 = \min\{i : f^{-n_i}(z) \notin S_1 \cup S_2\}$, and let $j_0 = \min\{j : m_j \geq n_{i_0}\}$. As before we can conclude that
$$\sum_{j=j_0}^{l} \ell(f^{-m_j}(J_\psi(z))) \leq C_2 C_3.$$

Fix some $z_1 \in S_1 \cap \Lambda$ and $z_2 \in S_2 \cap \Lambda$. Take $i < i_0$. Then $f^{-n_i}(z) \in S_1 \cup S_2$. Now assume that it is in $S_1$. Then, for every $m_j$ such that $n_i \leq m_j < n_{i+1}$, consider the box $B_{m_j} = f^{-(m_j-n_i)}(B(n_i))$, and $J_{m_j} = B_{m_j} \cap f^{-(m_j-n_i)}(J(z_1))$.

We have that $B_{m_j}$ has distortion $C_2$ and $B_{m_j} \cap B_{m_k} = \emptyset$ for every $0 \leq m_j, m_k < n_{i_0}$. Thus $J_{m_j} \cap J_{m_k} = \emptyset$. Therefore
$$\sum_{j=0}^{j_0-1} \ell(f^{-m_j}(J_\psi(z))) \leq \sum_{j=0}^{j_0-1} C_2 \ell(J_{m_j}) \leq 2C_2 K_1$$

where $K_1$ is as in Lemma 3.7.1. Set $C_4 = C_2 C_3 + 2C_2 K_1$. Then,
$$\sum_{j=0}^{l} \ell(f^{-m_j}(J_\psi(z))) \leq C_4,$$

and the claim is proved.

Finally, the proof of the lemma follows by the same arguments as in the proof of Lemma 3.5.2. □

We shall divide the proof of the Main Lemma into two cases: one, when $\Lambda$ is not a minimal set, and the other when it is. Remember that a compact invariant set is minimal if it has no properly compact invariant subset, or equivalently, if any orbit is dense.

3.7.1. *Case*: $\Lambda$ *is not a minimal set*.

LEMMA 3.7.4. *Let $B_\varepsilon(J)$ be a well-adapted box such that $\#\mathcal{R}(B_\varepsilon(J), \Lambda) = \infty$. Then there exists a return $\psi_0 \in \mathcal{R}(B_\varepsilon(J), \Lambda)$ such that the adapted box $B_{\psi_0} = \text{Image}(\psi_0)$ satisfies the conditions of Lemma 3.7.2; i.e., for every $\psi \in \mathcal{R}(B_{\psi_0}, \Lambda)$, $|\psi'| < \frac{1}{2}$ holds.*

*Proof.* Let $B_\varepsilon(J)$ be a well adapted box as in the hypothesis of the lemma, and let $K_1, K, C_1$ be as in Lemmas 3.7.1, 3.7.3 and Corollary 3.8 respectively. Consider also $L = \min\{\ell(J(z)) : z \in B_\varepsilon(J) \cap \Lambda\}$.

Let $r > 0$ be such that
$$r \frac{C_1}{L} \exp(K_0 K_1 + K_0 K) < \frac{1}{2}.$$

Since $\#\mathcal{R}(B_\varepsilon(J), \Lambda) = \infty$, there exists $\psi_0 \in \mathcal{R}(B_\varepsilon(J), \Lambda)$ such that

$$\ell(f^{-j}(J_{\psi_0}(z))) < r, \text{ for all } j \geq 0, \text{ for all } z \in B_{\psi_0} \cap \Lambda.$$

Let $k_0$ be such that $\psi_0 = f^{-k_0}_{/S_0}$, where $S_0 = \text{dom}(\psi_0)$.

Let us prove that the box $B_{\psi_0}$ satisfies the thesis of the lemma. Observe that if $z \in S_0 \cap \Lambda$, then for $y \in J(z)$

$$\|Df^{-k_0}_{/\tilde{F}(y)}\| \leq \frac{\ell(f^{-k_0}(J(z)))}{J(z)} \exp(K_0 K_1).$$

Let now $\psi \in \mathcal{R}(B_{\psi_0}, \Lambda), \psi = f^{-k}_{/S_\psi}, S_\psi = \text{dom}(\psi)$. Setting $n_0 = k - k_0$ ($k \geq k_0$), we have $f^{-n_0}(S_\psi) \subset S_0$.

Then, for $y \in J_{\psi_0}(z)$, $z \in \text{dom}(\psi)$,

$$\begin{aligned}
|\psi'(y)| &= \|Df^{-k}_{/\tilde{F}(y)}\| \leq \|Df^{-k_0}_{/\tilde{F}(f^{-n_0}(y))}\| \|Df^{-n_0}_{/\tilde{F}(y)}\| \\
&\leq \frac{\ell(f^{-k_0}(J(f^{-n_0}(z))))}{\ell(J(f^{-n_0}(z)))} \exp(K_0 K_1) \frac{\ell(f^{-n_0}(J_{\psi_0}(z)))}{\ell(J_{\psi_0}(z))} \exp(K_0 K) \\
&= \ell(f^{-n_0}(J_{\psi_0}(z))) \frac{\ell(J_{\psi_0}(f^{-k}(z)))}{\ell(J_{\psi_0}(z))} \frac{1}{\ell(J(f^{-n_0}(z)))} \exp(K_0 K_1 + K_0 K) \\
&\leq rC_1 \frac{1}{L} \exp(K_0 K_1 + K_0 K) < \frac{1}{2}.
\end{aligned}$$

The proof is finished. □

Now, we can prove the Main Lemma when $\Lambda$ is not a minimal set. We have to show that

$$\|Df^{-n}_{/F(z)}\| \to_{n \to \infty} 0$$

for every point $z \in \Lambda$.

Since $\Lambda$ is not a minimal set, there exists a point $x \in \Lambda$ such that $x \notin \omega(x)$. Take now a small well-adapted box $B_\varepsilon(J)$ associated to $x$ such that $B_\varepsilon(J) \cap \{f^n(x) : n \geq 1\} = \emptyset$. Then, since $\Lambda$ is transitive, we conclude that for this adapted box,

$$\#\mathcal{R}(B_\varepsilon(J), \Lambda) = \infty$$

(notice that if the point $x$ is a boundary point of $\tilde{\Lambda}$ the same conclusion holds). Now, by the previous lemma, there exists an adapted box $B_0$ such that satisfies the conditions of Lemma 3.7.2. This implies in particular that for every $y \in B_0 \cap \Lambda$,

$$\|Df^{-n}_{/F(y)}\| \to_{n \to \infty} 0.$$

Let $z$ be any point in $\Lambda$. There are two possibilities:

- The $\alpha$- limit set $\alpha(z)$ is properly contained in $\Lambda$. Then, $\alpha(z)$ is a hyperbolic set; thus
$$\|Df^{-n}_{/F(z)}\| \to_{n\to\infty} 0.$$

- $\alpha(z) = \Lambda$. Then, there exists $m_0$ such that $f^{-m_0}(z) \in B_0$, implying that
$$\|Df^{-n}_{/F(f^{-m_0}(z))}\| \to_{n\to\infty} 0$$
and so
$$\|Df^{-n}_{/F(z)}\| \to_{n\to\infty} 0.$$

This completes the proof of the Main Lemma in case $\Lambda$ is not a minimal set.

3.7.2. *Case: $\Lambda$ is a minimal set.* We begin by remarking that we cannot expect to use the same argument here as in the preceding case, due to the fact that if $\Lambda$ is a minimal set, then for every adapted box, the set of returns of this box is always finite. Nevertheless we shall exploit the fact that when $\Lambda$ is a minimal set there exist "boundary points."

LEMMA 3.7.5. *Assume $\Lambda$ is minimal set. Then, there exists an arbitrarily small adapted box $B_\varepsilon(J)$ such that $B_\varepsilon(J^+) \cap \Lambda = \emptyset$ or $B_\varepsilon(J^-) \cap \Lambda = \emptyset$.*

*Proof.* From Lemma 3.3.2 we get that, since $\Lambda$ is a minimal set, the center unstable and stable manifolds are dynamically defined.

Let $B_\varepsilon(J)$ be any small adapted box. We claim that
$$\bigcup_{z \in B_\varepsilon(J) \cap \Lambda} J(z) \neq B_\varepsilon(J)$$
and this happens on both sides of the box.

Arguing by contradiction, assume that
$$B^+_\varepsilon(J) = \bigcup_{z \in B_\varepsilon(J) \cap \Lambda} J(z).$$

Since $\Lambda$ is a minimal set, there exists some $n > 0$ sufficiently large such that $f^n(W^{cs}_\varepsilon(x)) \subset B^+_\varepsilon(J)$ $(x \in J)$. We can define a map
$$P : W^{cs,+}_\varepsilon(x) \to W^{cs,+}_\varepsilon(x)$$
in the following way: if $y \in W^{cs,+}_\varepsilon(x)$, then $f^n(y) \in f^n(W^{cs,+}_\varepsilon(x))$. Now there exists a point $z \in B^+_\varepsilon(J) \cap \Lambda$ such that $f^n(y) \in J(z)$. Let
$$P(y) = W^{cs,+}_\varepsilon(x) \cap J(z).$$

It is not difficult to see that the map $P$ is continuous, and we conclude that there exists a fixed point $y_0$ for this map. This means that there exists a point $z_0 \in B_\varepsilon^+(J) \cap \Lambda$ such that $y_0 \in J(z_0)$ and $f^n(y_0) \in J(z_0)$. Since $J(z_0) \subset W_\gamma^{cu}(z_0)$ and this central unstable manifold is dynamically defined, we conclude that $\alpha(z_0)$ is a periodic orbit, which is a contradiction to the fact that $\Lambda$ is a minimal set, which completes the proof of our claim.

Thus, as in Lemma 3.6.3, case 1, we can construct a well-adapted box $B_{\varepsilon'}(J)$ such that $B_{\varepsilon'}(J) = S_1 \cup B_{\varepsilon''}(J) \cup S_2$ with $S_i \cap \Lambda = \emptyset$.

With the same reasoning, it can be proved that

$$B_{\varepsilon'}(J) \neq \bigcup_{z \in B_{\varepsilon'}(J) \cap \Lambda} \left( W_\gamma^{cs}(z) \cap B_{\varepsilon'}(J) \right).$$

Then, there exists $x_0 \in B_{\varepsilon'}(J) \cap \Lambda$ such that there is no point of $\Lambda$ on one side of $W_\varepsilon^{cs}(x_0)$.

Finally, every small adapted box to this point $x_0$ satisfies the conclusion of the lemma. □

LEMMA 3.7.6. *Let $B_\varepsilon(J)$ be a well-adapted box such that $B_\varepsilon(J^+) \cap \Lambda = \emptyset$. Then $B_\varepsilon(J^+)$ is an "adapted box"; i.e., for all $y \in B_\varepsilon(J) \cap \Lambda$,*

$$f^{-n}(J^+(y)) \cap B_\varepsilon(J^+) = \emptyset \text{ or } f^{-n}(J^+(y)) \subset B_\varepsilon(J^+)$$

*where $J^+(y) = J(y) \cap B_\varepsilon(J^+)$. Moreover, there exists $K_1$ such that if $y \in B_\varepsilon(J) \cap \Lambda$ and $f^{-j}(J^+(y)) \cap B_\varepsilon(J^+) = \emptyset$, $1 \leq j < n$, then*

$$\sum_{j=0}^{n} \ell(f^{-j}(J^+(y))) < K_1.$$

*Proof.* Assume that, for some $y \in B_\varepsilon(J) \cap \Lambda$ and $n > 0$, $f^{-n}(J^+(y)) \cap B_\varepsilon(J^+) \neq \emptyset$ holds. As $B_\varepsilon(J)$ is an adapted box we conclude that $f^{-n}(J(y)) \subset B_\varepsilon(J)$. Moreover, there is a finite sequence of returns $\psi_i \in \mathcal{R}(B_\varepsilon(J), \Lambda)$, $i = 1, \ldots, k$, such that $f^{-n}(J(y)) = \psi_k \circ \cdots \circ \psi_1(J(y))$. If $f^{-n}(J^+(y))$ is not contained in $B_\varepsilon(J^+)$, then $f^{-n}(J^+(y)) \cap W_\varepsilon^{cs}(x) \neq \emptyset$. In particular, $W_\varepsilon^{cs}(x) \subset \text{Image}(\psi_k)$. Hence, $f^n(W_\varepsilon^{cs}(x) \cap J^+(y) \neq \emptyset$ and $f^n(W_\varepsilon^{cs}(x)) \subset \text{dom}(\psi_1) \subset B_\varepsilon(J)$. Therefore $f^n(W_\varepsilon^{cs}(x) \subset B_\varepsilon(J^+)$ and so does $f^n(x)$. Since $x \in \Lambda$ this is a contradiction, and completes the proof of the first part.

The existence of $K_1$ can be proved with the same arguments used in the proof of Lemma 3.7.3. □

LEMMA 3.7.7. *Let $B_\varepsilon(J)$ be a well-adapted box such that $B_\varepsilon(J^+) \cap \Lambda = \emptyset$. Then there exists $K$ such that for every $y \in B_\varepsilon(J) \cap \Lambda$,*

$$\sum_{j \geq 0} \ell(f^{-j}(J^+(y))) < K.$$

*In particular there exists* $J_1(y), J^+(y) \subset J_1(y) \subset J(y)$ *such that the length of* $J_1(y) - J^+(y)$ *is bounded away from zero* (*independently of* $y$) *and such that*

$$\sum_{n=0}^{\infty} \ell(f^{-n}(J_1(y))) < \infty.$$

*Proof.* First, we shall define returns of $B_\varepsilon(J^+)$ as we did for adapted boxes. Let $S$ be a vertical strip and define $S^+ = S \cap B_\varepsilon(J^+)$. A map $\psi : S^+ \to B_\varepsilon(J^+)$ is called a return of $B_\varepsilon(J^+)$ associated to $\Lambda$ if:

- $S \cap \Lambda \neq \emptyset$;
- there exists $k > 0$ such that $\psi = f^{-k}_{/S^+}$;
- $\psi(S^+) = f^{-k}(S^+)$ is a connected component of $f^{-k}(B_\varepsilon(J^+)) \cap B_\varepsilon(J^+)$;
- $f^{-i}(S^+) \cap B_\varepsilon(J^+) = \emptyset$ for $1 \leq i < k$.

If $\mathcal{R}(B_\varepsilon(J^+))$, the set of returns of $B_\varepsilon(J^+)$, is empty, by the preceding lemma, we conclude the proof. Assume that $\mathcal{R}(B_\varepsilon(J^+)) \neq \emptyset$, and let $\psi : S^+ \to B_\varepsilon(J^+)$ be a return, $\psi_{/S^+} = f^{-k}$. Let $B_\psi = \text{Image}(\psi)$.

We claim that $\text{Cl}(B_\psi) \cap W^{cs}_\varepsilon(x) = \emptyset$. Otherwise, the image of the segment $s = S \cap W^{cs}_\varepsilon(x)$ under $f^{-k}$ is contained in $W^{cs}_\varepsilon(x)$. This means in particular that there exists a point $y \in W^{cs}_\varepsilon(x)$ such that $f^k(y) \in W^{cs}_\varepsilon(x)$. Since $W^{cs}_\varepsilon(x)$ is in fact a stable manifold, we conclude that $\omega(x)$ is a periodic point, contradicting the fact that $\Lambda$ is a minimal set. This completes the proof of our claim.

If the $\#\mathcal{R}(B_\varepsilon(J^+))$ is finite, there exists $N$ such that for every $y \in B_\varepsilon(J) \cap \Lambda$ and $n \geq N$

$$f^{-n}(J(y)) \cap B_\varepsilon(J^+) = \emptyset.$$

To prove this, observe that, by the claim, there exists $r > 0$ such that

$$\text{dist}(W^{cs}_\varepsilon(x), \text{Image}(\psi)) > r, \text{ for all } \psi \in \mathcal{R}(B_\varepsilon(J^+)).$$

Since the central unstable manifolds are in fact unstable manifolds, there exists $N$ such that for every $y \in B_\varepsilon(J) \cap \Lambda$ and $n \geq N$ we have

$$\ell(f^{-n}(J(y))) < r.$$

Thus, if $f^{-n}(J(y)) \cap B_\varepsilon(J^+) \neq \emptyset$ for some $n \geq N$ then we get a contradiction. Indeed, if this happens, there exists a return $\psi \in \mathcal{R}(B_\varepsilon(J^+))$ such that $f^{-n}(J(y)) \cap \text{Image}(\psi) \neq \emptyset$, and so, by our definition of $r$, we conclude that $f^{-n}(J(y)) \subset B_\varepsilon(J^+)$ implying that $f^{-n}(y) \in B_\varepsilon(J^+) \cap \Lambda = \emptyset$, which is impossible.

Let $n_0 = \max\{n \geq 0 : f^{-n}(J(y)) \cap B_\varepsilon(J^+) \neq \emptyset\}$. Thus, $n_0 \leq N$ and, since $f^{-n_0}(J^+(y)) \subset J^+(f^{-n_0}(y))$, we conclude that

$$f^{-j}(J^+(f^{-n_0}(y))) \cap B_\varepsilon(J^+) = \emptyset, \text{ for all } j \geq 1$$

and so, by the previous lemma,
$$\sum_{j\geq 0} \ell(f^{-j}(J^+(f^{-n_0}(y)))) < K_1.$$

Therefore
$$\sum_{j\geq 0} \ell(f^{-j}(J^+(y))) < N\mathrm{diam}(M) + K_1 = K.$$

Finally, assuming that $\#\mathcal{R}(B_\varepsilon(J^+))$ is infinite, take $r > 0$ such that
$$\frac{r}{L}\exp(K_0 K_1) < \frac{1}{2}$$
where $L = \min\{\ell(J^+(y)) : y \in B_\varepsilon(J) \cap \Lambda\}$. Applying the same argument as before, we can conclude that there exists only a finite number of returns $\psi_1, \psi_2, \ldots, \psi_l$ such that $\mathrm{Image}(\psi_i)$ is not contained in a box $B_\varepsilon(J')$ of $x$ with "vertical size" at most $r$. In particular for every other return $\psi$, $\mathrm{Image}(\psi) \subset B_\varepsilon(J')$,
$$|\psi'| \leq \frac{r}{L}\exp(K_0 K_1) < \frac{1}{2}.$$

Thus, as in the proof of Lemma 3.7.2, there exists $\tilde{K}$ such that if for some $z \in B_\varepsilon(J) \cap \Lambda$ we have that
$$f^{-n}(J^+(z)) \cap \mathrm{dom}(\psi_i) = \emptyset, \ i = 1, \ldots, l,$$
then
$$\sum_{j\geq 0} \ell(f^{-j}(J^+(z))) \leq \tilde{K}.$$

On the other hand, as above, there exists $N$ such that for every $y \in B_\varepsilon(J) \cap \Lambda$ and $n \geq N$
$$f^{-n}(J^+(y)) \cap \mathrm{dom}(\psi_i) = \emptyset, \ i = 1, \ldots, l.$$

Take any $y \in B_\varepsilon(J) \cap \Lambda$. Consider, if it exists, $m_0 = \min\{n \geq N : f^{-n}(J^+(y)) \cap B_\varepsilon(J^+) \neq \emptyset\}$ and let $z = f^{-m_0}(y)$.

Hence,
$$\sum_{j\geq 0} \ell(f^{-j}(J^+(y))) \leq N\mathrm{diam}(M) + K_1 + \tilde{K} = K < \infty.$$

In particular, as in Schwartz's proof of the Denjoy theorem, we conclude that for all $y \in B_\varepsilon(J) \cap \Lambda$ there exist $J_1(y), J^+(y) \subset J_1(y) \subset J(y)$ such that the length of $J_1(y) - J^+(y)$ is bounded away from zero (independently of $y$) and such that
$$\sum_{n=0}^{\infty} \ell(f^{-n}(J_1(y))) < \infty$$

and the proof of the lemma is finished. □

Now, we can prove the Main Lemma when $\Lambda$ is a minimal set. We shall proceed as in the proof of Lemma 3.5.2. Using the notation of the preceding lemma, take
$$B = \bigcup_{y \in B_\varepsilon(J) \cap \Lambda} J_1(y).$$

Notice that $B$ is an open set of $\Lambda$, and for every $y \in B \cap \Lambda$ (i.e. $y \in J_1(y)$), we have
$$\sum_{n=0}^\infty \ell(f^{-n}(J_1(y))) < \infty$$

and so
$$\|Df^{-n}_{/F(y)}\| \to_{n\to\infty} 0.$$

Let $z$ be any point in $\Lambda$. Since $\Lambda$ is a minimal set there exist $m_0 = m_0(z)$ such that $f^{-m_0}(z) \in B$ and so
$$\|Df^{-n}_{/F(f^{-m_0}(z))}\| \to_{n\to\infty} 0,$$

implying that
$$\|Df^{-n}_{/F(z)}\| \to_{n\to\infty} 0.$$

This completes the proof of the Main Lemma. □

Universidad Federal do Rio de Janeiro, Rio de Janeiro, Brazil
*E-mail address*: enrique@impa.br

IMPA, Rio de Janeiro, Brazil
*E-mail address*: samba@impa.br